\documentclass[12pt,letterpaper,reqno]{amsart}
\usepackage{mathrsfs} 
\usepackage{amsmath}
\usepackage{amssymb}
\usepackage{amsthm}
\usepackage{indentfirst}
\usepackage{xspace}
\usepackage{multirow}
\usepackage{hyperref}
\usepackage{xcolor}
\usepackage{calc}
\usepackage{float}
\usepackage{makecell}
\definecolor{darkblue}{rgb}{0,0,0.4}
\hypersetup{pdfauthor={Emily Dumas and Andrew
    Neitzke},pdftitle={Opers and nonabelian Hodge: numerical studies},colorlinks=true,urlcolor=darkblue,linkcolor=darkblue,citecolor=darkblue}

\definecolor{mpl1}{HTML}{1F77B4}
\definecolor{mpl2}{HTML}{FF7F0E}

\usepackage{verbatim}
\usepackage[letterpaper,margin=25mm,marginparwidth=23mm]{geometry}
\usepackage{thmtools}
\usepackage[all]{xy}
\usepackage{capt-of}
\usepackage{enumitem}
\usepackage{graphicx}
\usepackage{siunitx}

\usepackage{booktabs}

\usepackage{placeins}
\usepackage{algorithm}
\usepackage{algorithmic}

\usepackage{xcolor}

\graphicspath{{figures/}}

\usepackage{marginnote}

\newlength{\floatpagetop}

\newlength{\bibitemsep}
\newlength{\bibparskip}
\let\oldthebibliography\thebibliography
\renewcommand\thebibliography[1]{%
  \oldthebibliography{#1}%
  \setlength{\parskip}{\bibitemsep}%
  \setlength{\itemsep}{\bibparskip}%
}

\setlist{itemsep = 0.25em, topsep = 0.25em}

\declaretheoremstyle[spaceabove=0.25cm,spacebelow=0.25cm,notefont=\normalfont\bfseries, notebraces={(}{)}]{theorem}
\declaretheoremstyle[spaceabove=0.25cm,spacebelow=0.25cm,bodyfont=\normalfont,notefont=\normalfont\bfseries, notebraces={(}{)}]{noital}
\declaretheoremstyle[spaceabove=0.25cm,spacebelow=0.25cm,bodyfont=\normalfont\color{darkgreen},notefont=\normalfont\bfseries, notebraces={(}{)}]{green}
\declaretheoremstyle[spaceabove=0.25cm,spacebelow=0.25cm,bodyfont=\normalfont,notefont=\normalfont\bfseries,qed=$\qedsymbol$,notebraces={(}{)}]{proofstyle}

\numberwithin{equation}{section}

\makeatletter
\def\paragraph{\smallskip\@startsection{paragraph}{4}%
  \z@\z@{-\fontdimen2\font}%
  {\normalfont\itshape}}
\makeatother

\renewcommand{\leq}{\leqslant}
\renewcommand{\geq}{\geqslant}

\usepackage{stackengine}
\usepackage{scalerel}

\begin{document}


\setcounter{page}{1}
\title[Parameter Slices of Rational Maps]{Emergence of Mandelbrot-like and Julia-like Structures in Parameter Slices of Rational Maps}
\author[Pedro Suarez]{Pedro Iván Suárez Navarro}

\date{\today}

 {\begin{abstract}
We study complex one-dimensional parameter slices in a three-parameter family of rational maps with two free critical points, obtained by imposing the existence of periodic orbits with prescribed multipliers. Using explicit parametrizations, we explore these slices numerically by analyzing the behavior of the critical orbits and approximating the corresponding connectedness loci. The computations reveal rich parameter space structures closely analogous to those arising in cubic polynomial families, including Mandelbrot-like sets. In addition, we observe regions exhibiting Julia-like structures embedded in parameter space, arising from the interaction between bounded and escaping critical orbits. 
While the appearance of such structures is well established in polynomial dynamics, it remains comparatively less explored in the setting of rational maps. Our results provide numerical evidence that these parameter slices contain subsets closely related to the period-one and period-two slices of cubic polynomial families.
More precisely, certain regions appear to exhibit geometric and dynamical features consistent with embedded copies of these classical parameter spaces. These observations highlight how classical structures from polynomial dynamics can emerge naturally within parameter slices of rational maps.
\end{abstract}}

\maketitle


\section{Introduction}
The appearance of Mandelbrot-like and Julia-like structures in parameter spaces is a recurring phenomenon in complex dynamics. In polynomial families, Mandelbrot-like structures are supported by a well-developed theoretical framework, notably through polynomial-like mappings, renormalization, and general results on bifurcation loci (see, e.g., \cite{mcmullen2000mandelbrot}).

In contrast, the occurrence of Julia-like structures within parameter spaces remains far less understood. While such phenomena have been rigorously established in specific settings, most notably in cubic polynomial families and in certain constructions involving rational maps, their presence in more general contexts is still only partially explored.

Beyond these classical cases, similar patterns are frequently observed in numerical experiments across a wide range of dynamical systems. Their ubiquity is often regarded as part of the folklore of the field: they arise consistently in computations and specific examples, yet lack a unified theoretical description.

This motivates the question of whether such structures persist in broader classes of dynamical systems, and to what extent they can be detected, organized, and interpreted through concrete computational and geometric approaches, as illustrated by the numerical experiments presented in this work.

In this paper, we present numerical experiments on complex one-dimensional parameter slices arising in a family of rational maps with two superattracting fixed points and a periodic point with prescribed complex multiplier. A standard approach to the study of parameter spaces is through parameter slices, obtained by imposing dynamical constraints on periodic or critical orbits. In particular, Milnor introduced curves of the form $\mathrm{Per}_n(\lambda)$, consisting of maps that possess a periodic orbit of period $n$ with multiplier $\lambda$ \cite{milnor1993geometry,milnor2000rational}. These curves play a central role in organizing the global geometry of parameter spaces and have been extensively studied in the context of cubic polynomials with a periodic critical orbit \cite{milnor1992remarks,milnor20099,bonifant2010cubic,bonifant2025cubic}.

\begin{figure}[ht!]
\centering
\includegraphics[width=0.9\textwidth, height=0.5\textheight]{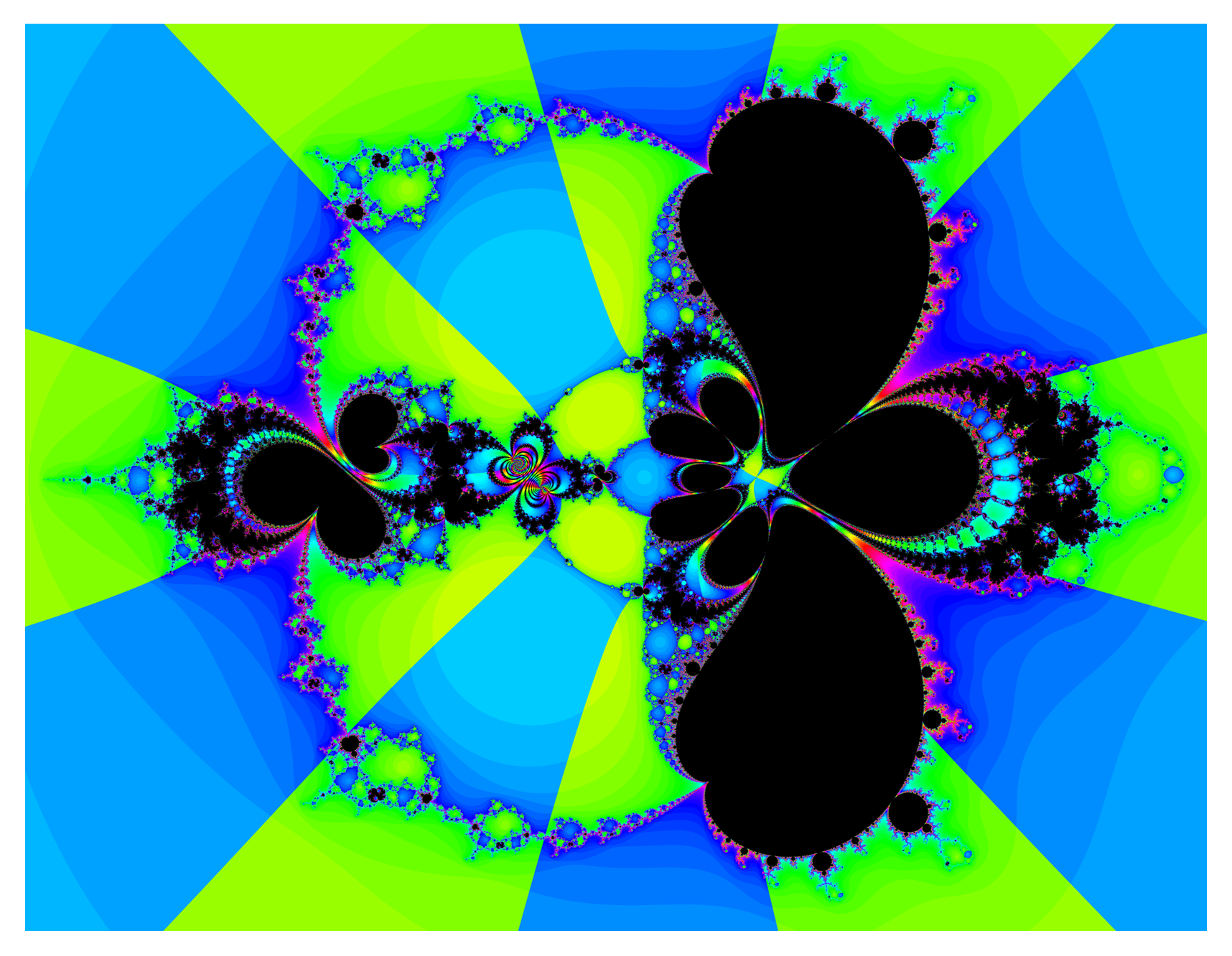} 
\hspace{0.02\textwidth}
 \includegraphics[width=0.44\textwidth, height=0.3\textheight]{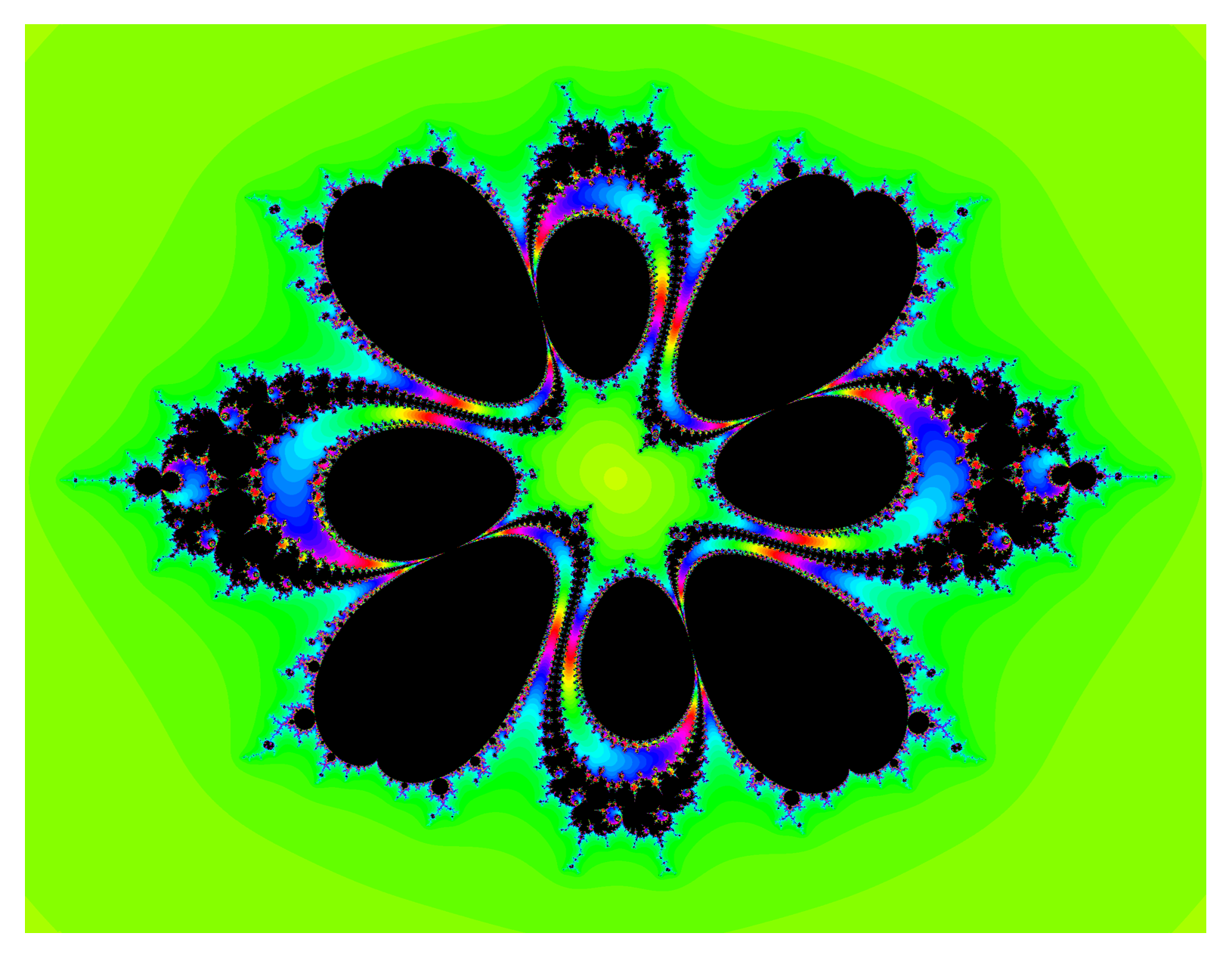}
\hspace{0.02\textwidth}
   \includegraphics[width=0.44\textwidth, height=0.3\textheight]{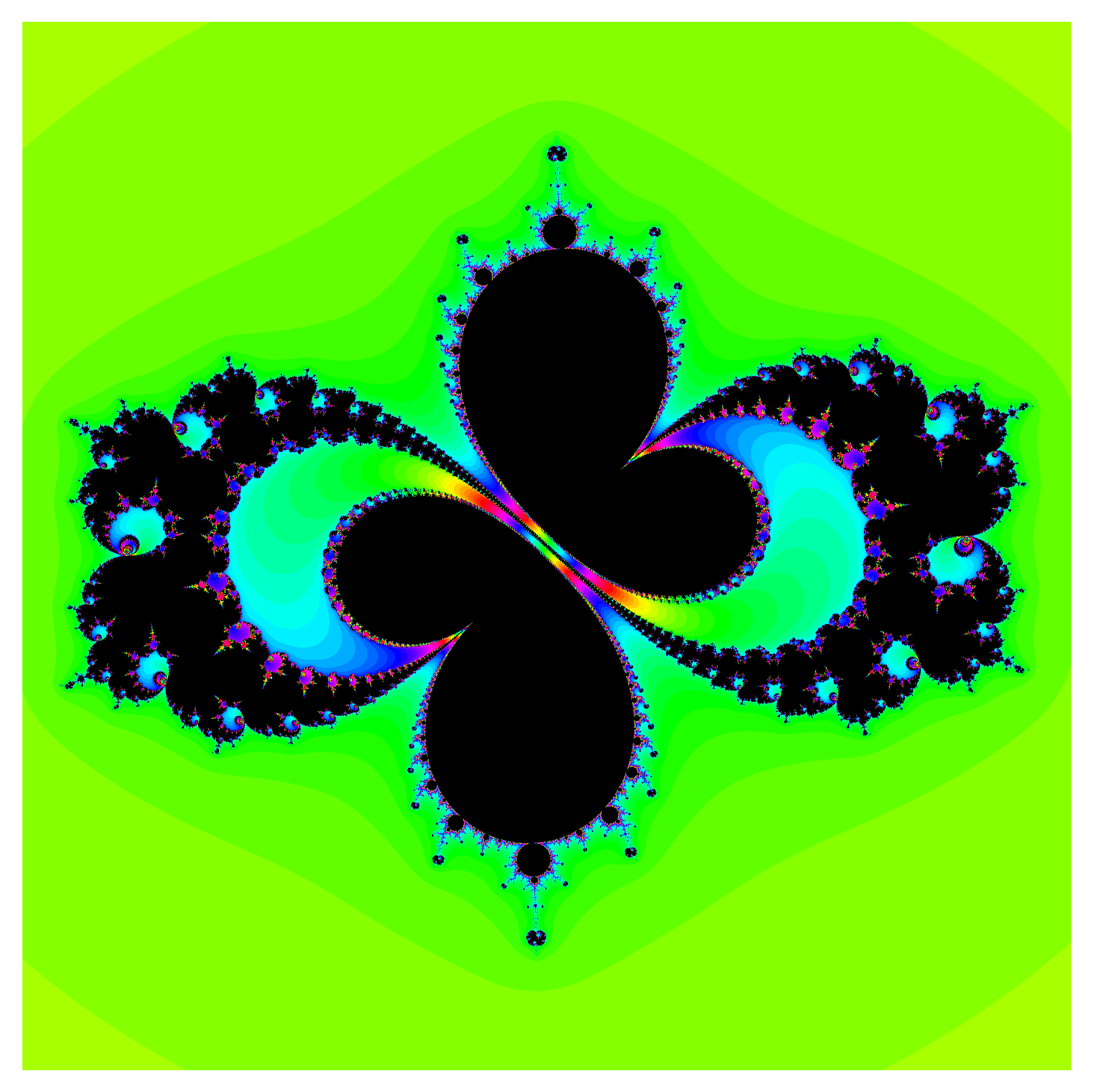}
\caption{
Comparison between parameter slices in the rational family studied in this paper (top) and in the cubic polynomial family (bottom). The close visual resemblance, including Mandelbrot-like and Julia-like structures, suggests that parts of the rational parameter slices may exhibit dynamics analogous to those of cubic polynomial families.
}
\label{fig:intro_comparison}
\end{figure}

Our work is motivated by a conjectural picture predicting the presence of embedded copies of parameter slices of cubic polynomials
\[
P_{\mu,b}(z) = \mu z + b z^2 + z^3, \quad (\mu, b) \in \mathbb{C}^2,
\]
within the connectedness loci of the rational families considered here. Related phenomena arise in the neutral case $\mu = e^{2\pi i p/q}$ \cite{zakeri1999dynamics, zhang2025parabolic, zhang2024dynamical}. In particular, Buff and Henriksen showed that, for $|\lambda|=1$, the bifurcation locus of the corresponding one-parameter family in $b$ contains quasiconformal copies of $J(\lambda z + z^2)$~\cite{buff2001julia}, a result later extended to higher-degree polynomials in~\cite{inou2012combinatorics}. Moreover, even in the classical Mandelbrot set of the quadratic family, quasiconformal copies of Julia sets (in particular, Cantor Julia sets) are known to appear~\cite{KawahiraKisaka2026}.

These phenomena are closely related to the theory of polynomial-like mappings \cite{DouadyHubbard} and renormalization, which provide a natural framework for understanding the emergence of Mandelbrot-like structures in parameter space. In the parabolic setting, analogous structures are expected to be described by the theory of parabolic-like mappings \cite{lomonaco2015parabolic}, reflecting the role of parabolic renormalization in organizing the geometry of parameter spaces near parabolic parameters.

The experiments reported here focus on periods one and two. These cases remain computationally tractable while capturing transitions between different regions of parameter space, enabling a systematic exploration of the dependence on the multiplier. The present work is also inspired by the visual structures observed in \cite{baranski1998newton, baranski2001newton} for cubic rational maps $z \mapsto a z^2\frac{b z+1-2b}{(2-b)z-1}$ with $a \in \mathbb{C}\setminus \lbrace  0\rbrace$ and $b \in \mathbb{C}\setminus \lbrace  0,1\rbrace$, as well as in \cite{buff2013limits} for parabolic quadratic rational maps, e.g. $z \mapsto \frac{\rho z}{1+a z+z^2}$ with $a \in \mathbb{C}$ and $\rho \in \mathbb{C}^{*} 
$.

The numerical explorations carried out in this work reveal recurring geometric patterns in the parameter slices, suggesting a strong connection with structures arising in cubic polynomial families.

\noindent\textbf{Conjectural framework.}
The numerical experiments suggest that the parameter slices $\mathscr{S}_n(\lambda)$, for $n=1,2$, exhibit geometric and dynamical structures closely analogous to those in slices of cubic polynomial parameter spaces. In particular, certain regions appear to contain subsets that are dynamically (or hybrid) equivalent to such spaces, arising as embedded copies within the connectedness loci of the rational family. These observations suggest that such phenomena may occur along loci where one critical orbit remains bounded while the other escapes, leading to configurations reminiscent of Julia sets in parameter space (cf. Figure~\ref{fig:intro_comparison}).

This raises several natural questions, including whether these loci are bounded and connected for given values of $\lambda$, whether such structures occur systematically within the slices $\mathscr{S}_n(\lambda)$, and whether subsets that appear to be homeomorphic copies of the Mandelbrot set (including baby Mandelbrot sets and parabolic Mandelbrot sets) are present. It is also natural to ask whether these observations can be made precise within a rigorous framework, and to what extent dynamically defined invariant sets in the dynamical plane may manifest themselves in parameter space through the interaction of critical orbits.

The goal of this paper is not to resolve these questions, but rather to provide numerical evidence supporting this picture and to illustrate the geometric structures that arise in these parameter slices. We also hope that these observations will invite the reader to reflect on possible theoretical directions suggested by these simulations, and to formulate more precise questions.

The article is organized as follows. In Section~\ref{sec:obst}, we introduce the family of rational maps under consideration, describe the corresponding parameter slices, and derive explicit parametrizations. Section~\ref{sec:plts} contains the numerical experiments and illustrates the geometric structures suggested by the simulations.

\section{The Family and Its Parameter Slices}
\label{sec:obst}

\subsection{The family and critical points}

We consider a family of rational maps of degree $d+1$ with two superattracting fixed points and two simple free critical points. After conjugation by a Möbius transformation, we normalize the fixed points to $c_1 = 0$ and $c_2 = \infty$, obtaining the following family

\begin{equation}\label{maineq}
f(z)=z^d \frac{z+\alpha}{\beta z+\gamma}, \qquad 
\alpha,\beta,\gamma\in\mathbb{C}, \quad \gamma-\alpha\beta\neq 0.
\end{equation}

The map $f$ has exactly two free critical points given by
\[
z_{1,2} =
-\frac{(-1+d)\alpha\beta + (1+d)\gamma
\pm
\sqrt{(-\alpha\beta+\gamma)\big((1+d)^2\gamma-(d-1)^2\alpha\beta\big)}}{2d\beta}.
\]
The dynamics is largely governed by the behavior of these critical orbits, which play a central role in organizing the parameter space.

To investigate the global structure of the family, we impose dynamical conditions on periodic orbits and consider the corresponding parameter slices. For each $n\geq 1$ and $\lambda\in\mathbb{C}$, we define the parameter slices
$\mathscr{S}_n(\lambda)$
as the set of conjugacy classes of maps in the family possessing a periodic orbit of period $n$ with multiplier $\lambda$.
 These conditions reduce the effective dimension of the parameter space and produce complex one-dimensional slices that can be explored computationally.

In practice, the parameters $(\alpha,\beta,\gamma)$ can be expressed in terms of a single complex variable subject to the imposed dynamical constraints, allowing direct sampling of the parameter plane and a systematic investigation of the associated critical dynamics; in this way, these slices provide a natural setting for comparing the geometry of rational parameter spaces with that of cubic polynomial families.

The full parameter space is a subset of $\mathbb{C}^3$ and its global structure is highly intricate. Even in low-degree cases (see \cite{baranski1998newton}, Theorem~3.7, for $d=2$), the geometry already exhibits considerable complexity. Restricting to one--dimensional parameter slices provides a tractable framework for computation, allowing effective sampling, visualization, and a systematic study of the dynamical behavior.

Within these parameter slices, we focus on the subset of parameters for which the dynamics remains bounded. We are particularly interested in the connectedness locus $\mathcal{M}_{\lambda}$ of $\mathscr{S}_n(\lambda)$, defined by
\[
\mathcal{M}_{\lambda}
=
\left\{ f \in \mathscr{S}_n(\lambda) : \mathscr{J}(f) \text{ is connected} \right\}.
\]

Numerical experiments for different values of $\lambda$ with $\vert \lambda \vert \leq 1$ suggest that $\mathcal{M}_{\lambda}$ is bounded and connected.

To explore these sets computationally, we derive explicit parametrizations of the parameter slices introduced above. The goal is to reduce the parameter space to one complex dimension and to enable systematic exploration via the behavior of critical orbits.

\subsection{Parameter slices $\mathscr{S}_n(0)$}
\label{par_pern(0)}
Let $f$ be a map in the family defined in~\eqref{maineq}. If $\alpha = \beta = 0$, then $f$ is conjugate to the monomial $z^{d+1}$. Otherwise, the map has two free critical points. By a conformal change of coordinates, we may assume that one of them is located at $z = 1$, which imposes the relation
\[
\beta(d+(d-1)\alpha) + (1 + d(1+\alpha))\gamma = 0.
\]

Assuming $\gamma \neq 0$, we obtain
\[
\beta = -\frac{(1 + d(1+\alpha))\gamma}{d + (d-1)\alpha}.
\]

Setting
\[
a = -\frac{d + (d-1)\alpha}{\gamma}, 
\qquad 
b = \frac{1}{d + (d-1)\alpha},
\]
and substituting into~(2.1), we obtain the normalized family
\begin{equation}\label{baranskifamily}
    \begin{aligned}
        f_{a,b}(z) = a z^d 
\left(
\frac{bz + \frac{1-db}{d-1}}{\frac{d-b}{d-1}z - 1}
\right),
\qquad a \in \mathbb{C}^\times, b \in \mathbb{C}.
    \end{aligned}
\end{equation}

The map $f_{a,b}$ has two critical points
\begin{equation}\label{critPoints}
    \begin{aligned}
z_1 = 1, 
\qquad 
z_2 = \frac{1-db}{b(b-d)},
    \end{aligned}
\end{equation}
with corresponding critical values
\[
v_1 = a, 
\qquad 
v_2 = \frac{(1-bd)^{d+1}}{b^{d-1}(b-d)^{d+1}}\,a.
\]

This normalization reduces the study of $\mathscr{S}_n(0)$ to a one-dimensional parameter space. More precisely, for each $d \geq 2$, we obtain the following parametrizations:

\subsubsection*{Period 1}
The condition that $z = 1$ is a superattracting fixed point is equivalent to $a = 1$, leaving $b \in \mathbb{C}$ as a free parameter. Substituting into~\eqref{maineq} then yields a one-parameter family $f_b$.

\subsubsection*{Period 2}
The condition that $z = 1$ belongs to a superattracting $2$--cycle, that is, $f_{a,b}^{2}(1) = 1$, imposes
\begin{equation}\label{eqbb}
b=\frac{1-a^{d+1}+(a-1)^d}{a\bigl(1-a^{d+1}-d(1-a)a^d\bigr)}.
\end{equation}
Substituting into~\eqref{maineq}, we obtain a one-parameter family $f_a$.

The parametrization of $\mathscr{S}_2(0)$ admits a finite set of singular parameters associated with degeneracies. Numerically, these appear to act as organizing centers for nearby structures in the parameter plane.

\subsection{Parametrization of $\mathscr{S}_1(\lambda)$}
\label{subsec:lambda1}
We now fix $\lambda \in \mathbb{C}$ and consider parameter slices with a fixed point of multiplier $\lambda$. Without loss of generality, we assume that this fixed point is located at $z = 1$. Under this normalization, the parameters satisfy
\[
\gamma = 1 + \alpha - \beta.
\]
Imposing $f'(1) = \lambda$ yields
\[
1 + d + d\alpha - \beta - \lambda - \alpha \lambda = 0.
\]

Solving for $\alpha$ and $\gamma$, we obtain
\[
\alpha = \frac{-1 - d + \beta + \lambda}{d - \lambda}, 
\qquad 
\gamma = \frac{-1 + \beta - d\beta + \beta \lambda}{d - \lambda}.
\]

Substituting into~\eqref{maineq} yields a one-parameter family $f_{\beta,\lambda}$,
which provides a parametrization of the slice $\mathscr{S}_1(\lambda)$.

\subsection{Parametrization of $\mathscr{S}_2(\lambda)$}
\label{subsec:lambda2}
We now fix $\lambda \in \mathbb{C}$ and impose the existence of a $2$--cycle with multiplier $\lambda$. More precisely, we require that $z = 1$ belongs to a $2$--cycle, that is,
\[
f^{2}(1) = 1, \qquad (f^{2})'(1) = \lambda.
\]
Writing $t = f(1)$, these conditions become
\[
f(1) = t, \qquad f(t) = 1, \qquad f'(1)f'(t) = \lambda,
\]
which leads to the system
\[
\begin{cases}
1 + \alpha - t\beta - t\gamma = 0, \\
t^{1+d} + t^d \alpha - t\beta - \gamma = 0.
\end{cases}
\]
Substituting into the multiplier condition yields a quadratic equation in $\alpha$. Solving for $\alpha$ as a function of $t$ (and $\lambda$) and substituting into~\eqref{maineq}, we obtain a parametrization of the slice $\mathscr{S}_{2}(\lambda)$ in terms of $t$, leading to a family $f_{t,\lambda}$.

\section{Numerical experiments}
\label{sec:plts}

We explore the geometry of the parameter slices by systematically sampling the corresponding parameter planes and analyzing the dynamics of the free critical orbits, using the parametrizations derived in Section~\ref{sec:obst}.

The goal of these experiments is not only to visualize the parameter slices, but also to identify recurring geometric patterns and relate them to known structures in polynomial dynamics. In particular, we focus on the interaction between bounded and escaping critical orbits, which plays a central role in organizing the parameter space.

For comparison, we refer to analogous parameter slices in cubic polynomial families, typically described using standard normal forms and parametrizations in the literature (see \cite{Milnor2000,inouvisualization,zakeri1999dynamics,zhang2025parabolic}). These provide a natural point of comparison for interpreting the structures observed in our rational slices.

For each parameter in the different slices under consideration, we iterate both critical points and classify parameters according to whether the corresponding orbits remain bounded. The parameter plane is sampled on a uniform grid with step size $10^{-4}$, and computations are performed with a maximum of $N = 4000$ iterations, implemented in \textit{Mathematica}.

To visualize the resulting structures and approximate the connectedness locus, we assign colors according to the behavior of the critical orbits: escape to one of the superattracting fixed points $0$ or $\infty$, attraction to the marked critical point $z = 1$, or convergence to another attracting or parabolic cycle distinct from $z = 1$.

To detect bounded dynamics, we use a forward-invariant trapping region
\[
A = \left\{ z \in \mathbb{C} :
\min\left\{\frac{|\alpha|}{2}, \frac{|\gamma|}{2|\beta|}, \left(\frac{|\gamma|}{3|\alpha|}\right)^{\frac{1}{d-1}} \right\}
\le |z| \le
\max\left\{2|\alpha|, \frac{2|\gamma|}{|\beta|}, (3|\beta|)^{\frac{1}{d-1}} \right\}
\right\}.
\]
A parameter is classified as belonging to the connectedness locus if both critical orbits enter the region $A$ and remain bounded under iteration up to $N$ steps. Outside $A$, escape is detected using a standard escape-time criterion.

\subsection{Slices $\mathscr{S}_n(0)$}
\FloatBarrier

We begin by investigating the slices $\mathscr{S}_n(0)$ for different values of the degree $d$. Using the parametrizations derived in Subsection~\ref{par_pern(0)}, we explore the geometry of the corresponding connectedness loci and analyze how their structure evolves as $d$ varies, with particular attention to recurrent patterns and their relation to cubic polynomial parameter spaces.

We focus on the cases $n=1$ and $n=2$, where the parameterizations remain sufficiently explicit to allow for efficient numerical exploration, typically in rational form.

For higher periods, the situation becomes substantially more involved: already for $n=3$, obtaining a one-parameter description of $\mathscr{S}_n(0)$ requires eliminating variables from increasingly complicated algebraic relations, making both symbolic and numerical analysis difficult.

Alternative approaches, such as parameterizations arising from Hamiltonian structures in the spirit of ~\cite{bonifant2010cubic}, may provide a more tractable framework for higher-period slices; we leave this as an open direction.

\subsubsection{Slices $\mathscr{S}_1(0)$}
For low values of $d$, the connectedness locus of $\mathscr{S}_{1}(0)$ already exhibits a central connected region together with escape components associated with the two superattracting fixed points of the family. In the case $d=2$, this locus coincides with the parameter space of cubic rational maps arising from Newton's method.

\begin{figure}[ht!]
   \centering   
   \includegraphics[width=0.4\textwidth, height=0.3\textheight]{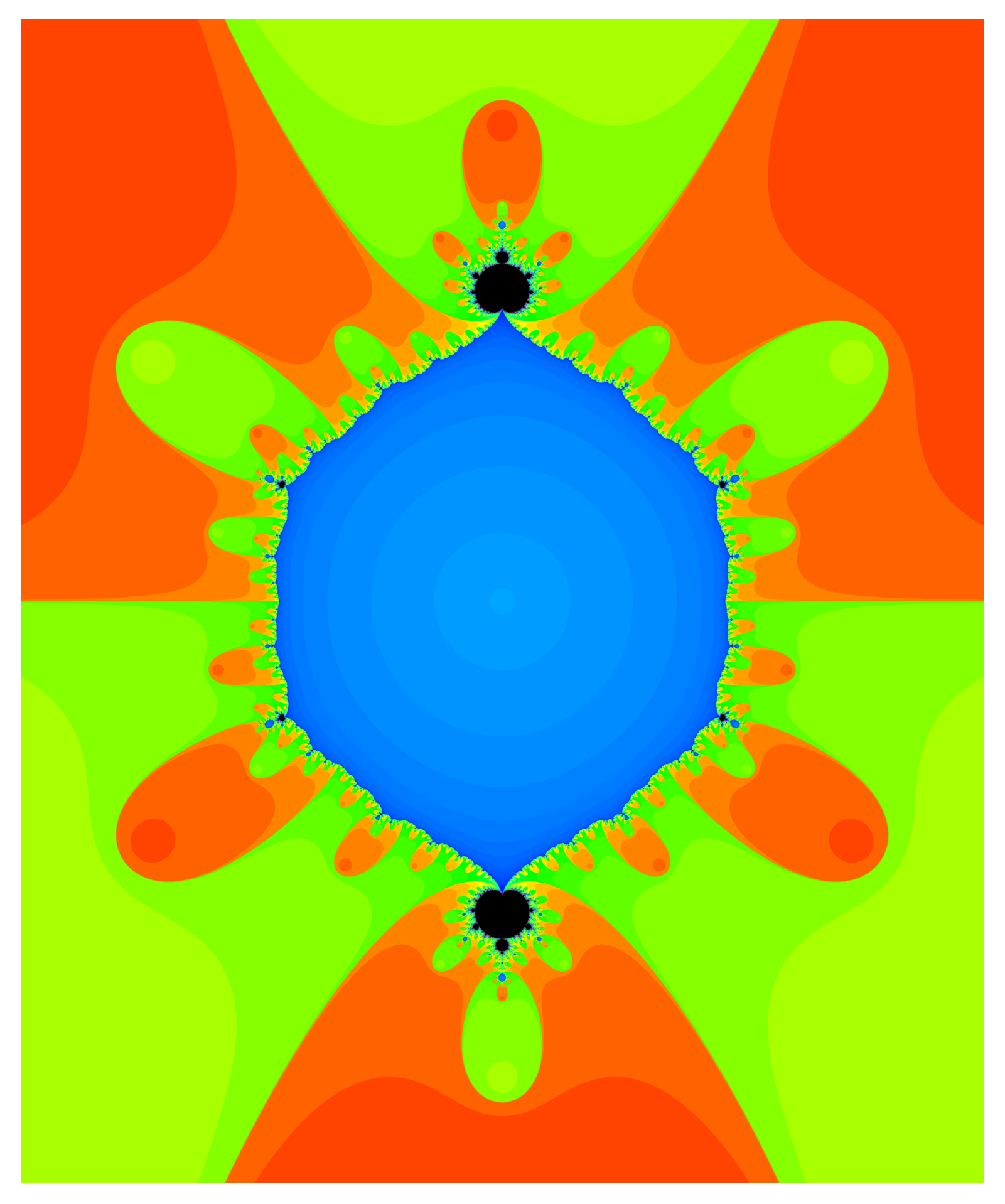} 
   \includegraphics[width=0.44\textwidth, height=0.3\textheight]{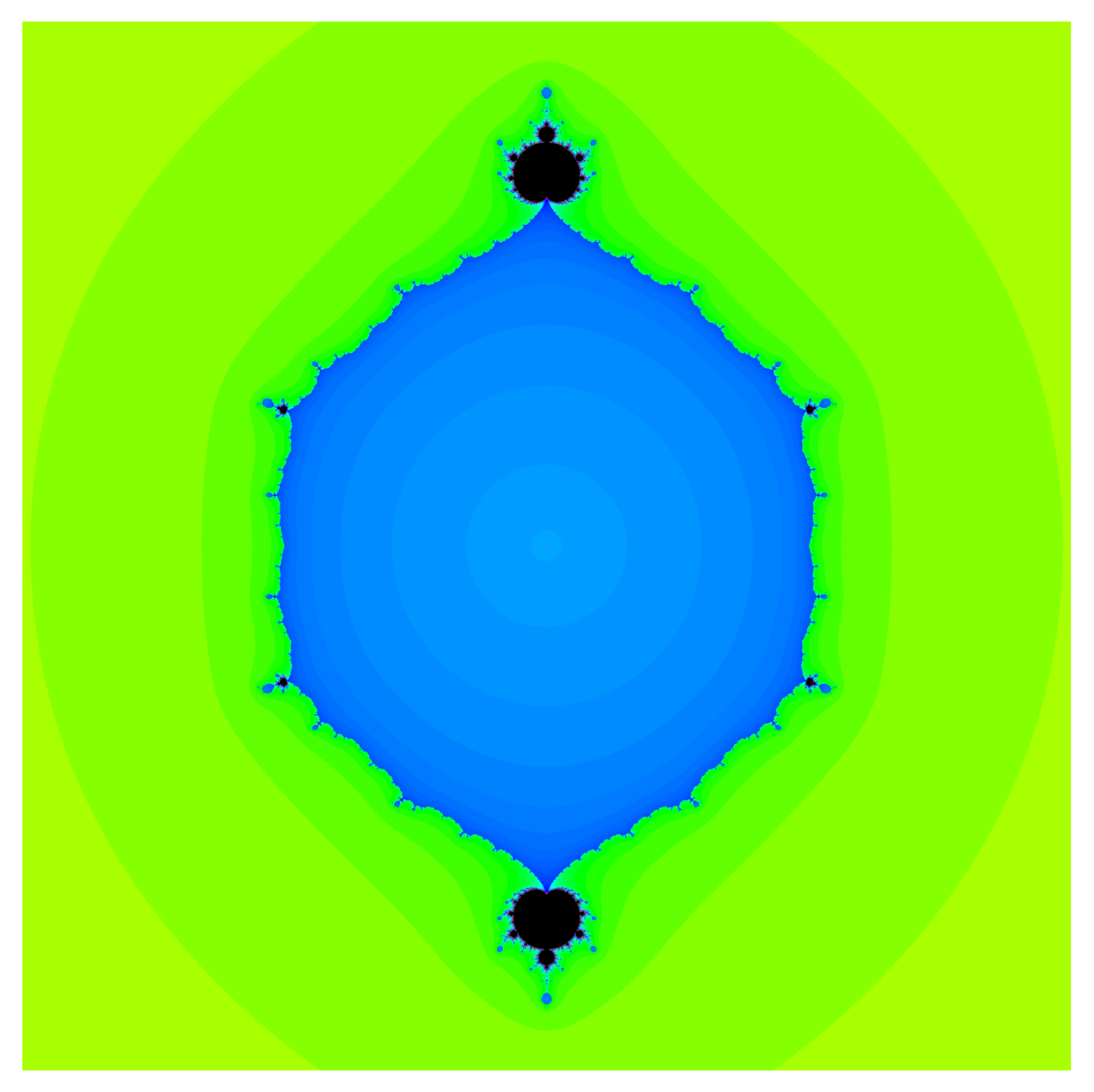} 
   \caption{Connectedness locus of $\mathscr{S}_1(0)$ for large $d$, compared with the cubic slice $\mathrm{Per}_1(0)$.}
   \label{fig:secper10Large}
\end{figure}

For $d \geq 3$, the central region takes on a lemon-like shape that becomes progressively flatter as $d$ increases. As shown in Figure~\ref{fig:secper10Large}, this region becomes increasingly similar to the slice $\mathrm{Per}_1(0)$ in the cubic polynomial family. Additional small components appear attached to the main body, clustering near singular parameters. This behavior is consistent with structures observed in cubic polynomial families (see, for instance, \cite{milnor20099,bonifant2023relations,petersen2025lemon}). The observed symmetry reflects the choice of parameter plane under the Möbius transformation $b \mapsto \frac{b - i}{b + i}$.

\subsubsection{Slices $\mathscr{S}_2(0)$}

We use the parametrization~\eqref{eqbb} and focus on degrees $d \geq 3$, where new geometric features arise. 
The case $d=2$ corresponds to the period-two slice in the parameter space of cubic Newton maps, which has been extensively studied (see \cite{baranski1998newton}) and will not be considered here.

\begin{figure}[ht!]
   \centering   
   \includegraphics[width=0.8\textwidth, height=0.5\textheight]{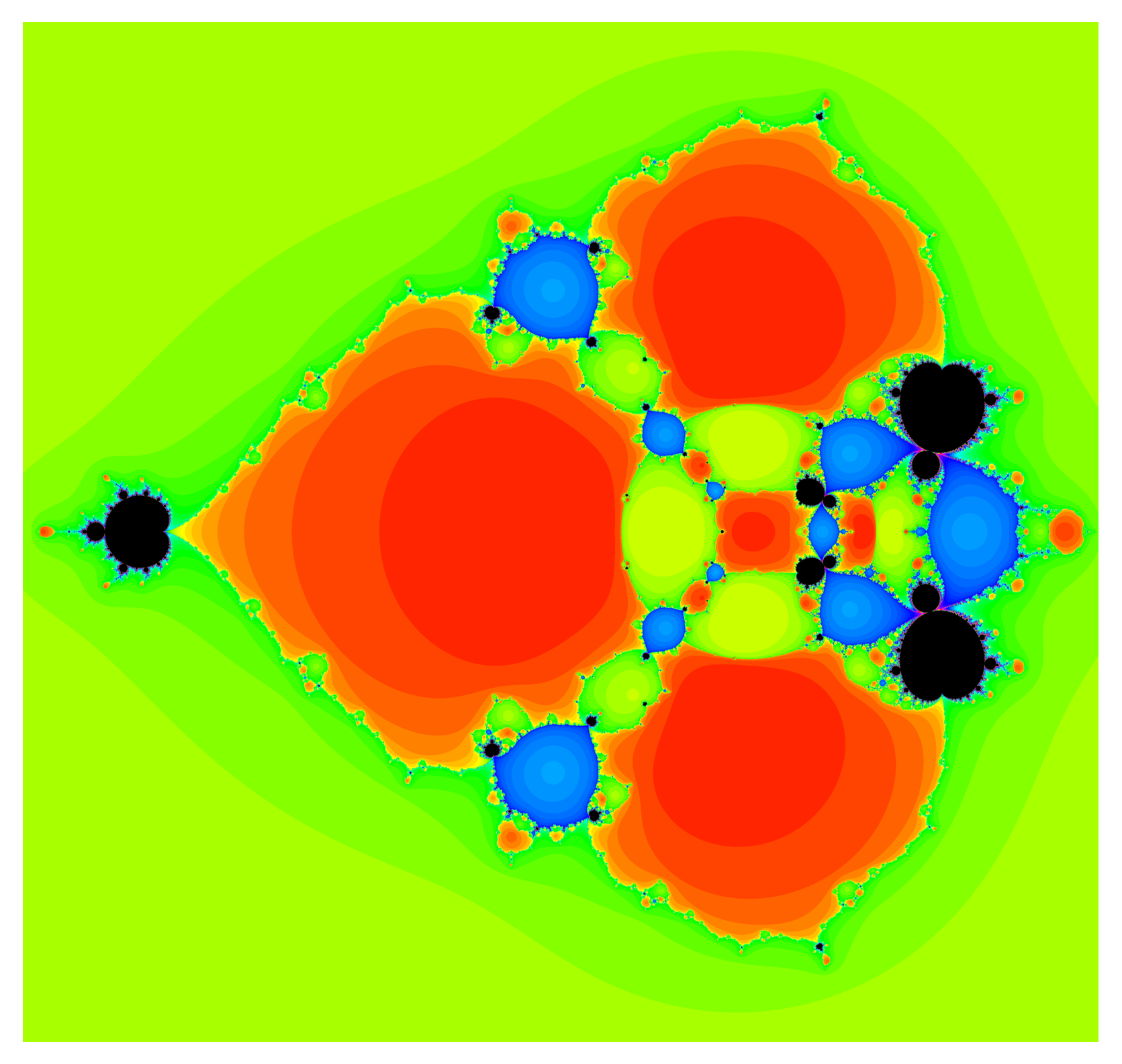} 
   \caption{Connectedness locus of $\mathscr{S}_2(0)$ in the $a$-parameter plane for $d=3$.} 
   \label{fig:per2_d3}
\end{figure}

Figure~\ref{fig:per2_d3} shows the connectedness locus of $\mathscr{S}_2(0)$ in the $a$-parameter plane for $d=3$. 
The picture reveals a decomposition into hyperbolic components, each containing a distinguished center corresponding to a singular parameter (puncture). 

These punctures arise from degeneracies of the parametrization \eqref{eqbb}, where the map fails to be well-defined or the dynamics becomes singular. 
In particular, they are determined by algebraic conditions on \eqref{eqbb}, such as $b=0$, $b=\infty$, or $b=1$. 

Dynamically, these correspond to configurations where the map loses degree. 
In addition, other singular centers associated with different components appear when the critical orbit of $z_2$ (see~\ref{critPoints}) collides with $0$, $\infty$, or $z_1$. 

For each $d$, this produces a finite collection of distinguished parameters, whose number grows linearly with the degree. 
Numerically, these points appear to organize the surrounding geometry: the visible hyperbolic components are arranged around them, and finer filamentary structures tend to accumulate near their locations. 
This suggests that the global structure of the slice is strongly influenced by these singular parameters and by the interaction between the critical orbits.

\begin{figure}[ht!]
   \centering   
   \includegraphics[width=0.61\textwidth, height=0.52\textheight]{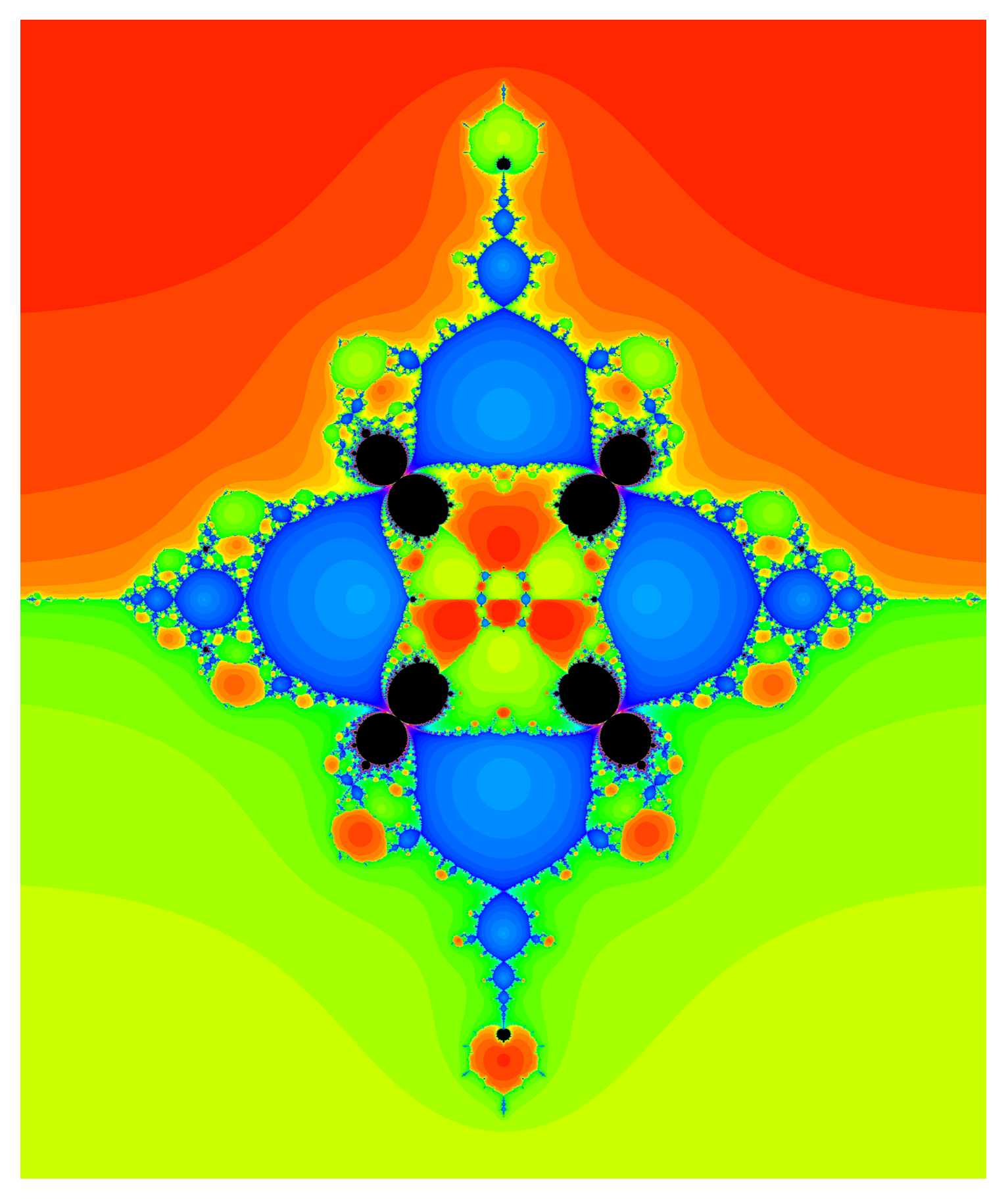}    
   \includegraphics[width=0.34\textwidth, height=0.52\textheight]{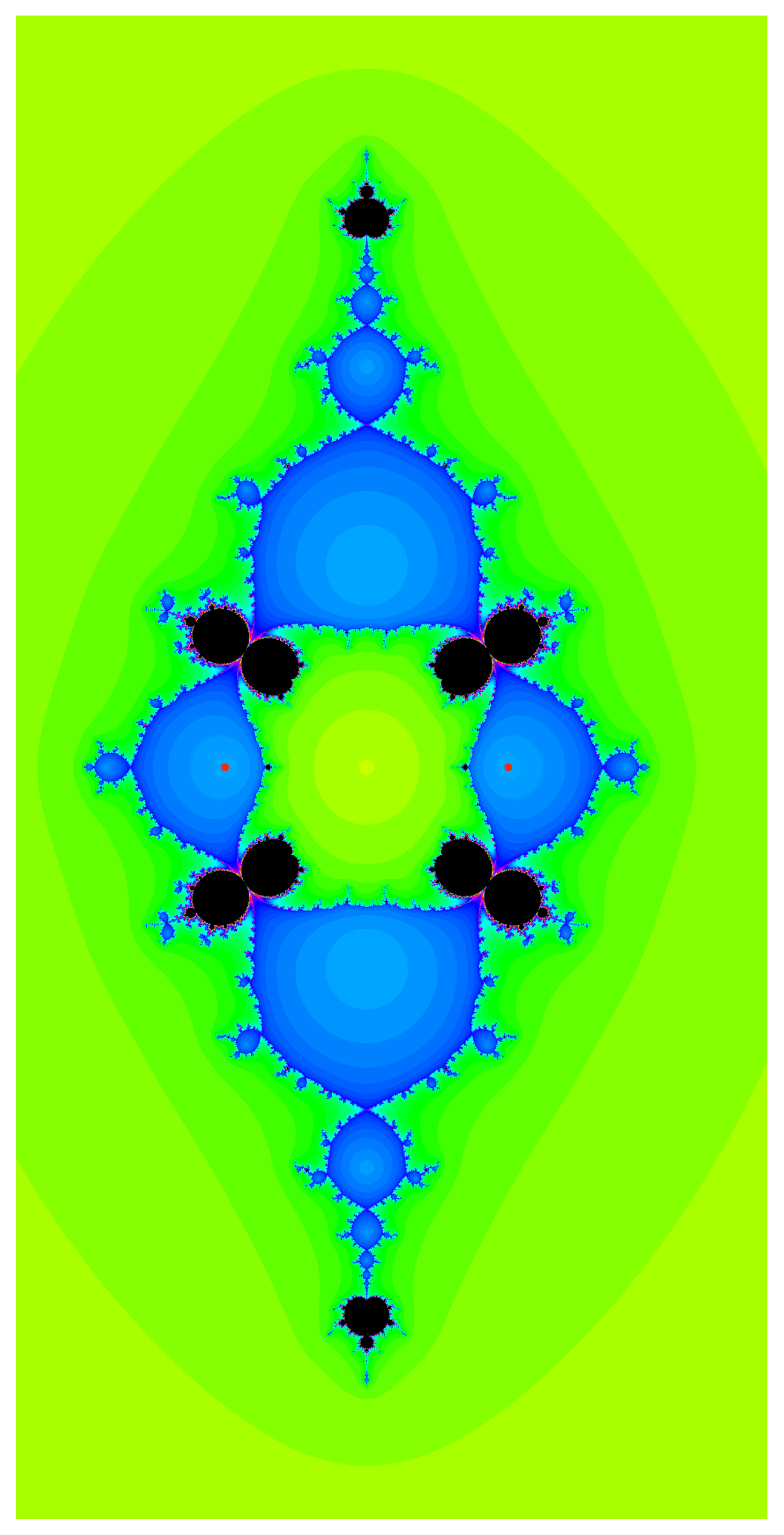} 
   \caption{Connectedness locus of $\mathscr{S}_2(0)$ for $d=3$ in the parameter plane given by $\frac{a - i}{a + i}$, illustrating its similarity with the corresponding slice $\mathrm{Per}_{2}(0)$ of cubic polynomials.} 
   \label{fig:secper20vertd3}
\end{figure}

Moreover, under the change of coordinates $a \mapsto \frac{a - i}{a + i}$, the similarity with the cubic slice becomes more apparent. 
Figure~\ref{fig:secper20vertd3} further illustrates this correspondence with the analogous cubic parameter space.

This phenomenon persists for larger values of $d$. 
Numerical experiments suggest a proliferation of components of type $A$ and $B$ as $d$ increases, while the overall organization near $a = 1$ remains remarkably stable. 
Moreover, the observed structures exhibit components consistent with all types in the sense of Milnor’s classification~\cite{milnor20099}, namely types $A$, $B$, $C$, and $D$.

\begin{figure}[ht!]
\centering
\includegraphics[width=0.35\textwidth, height=0.2\textheight]{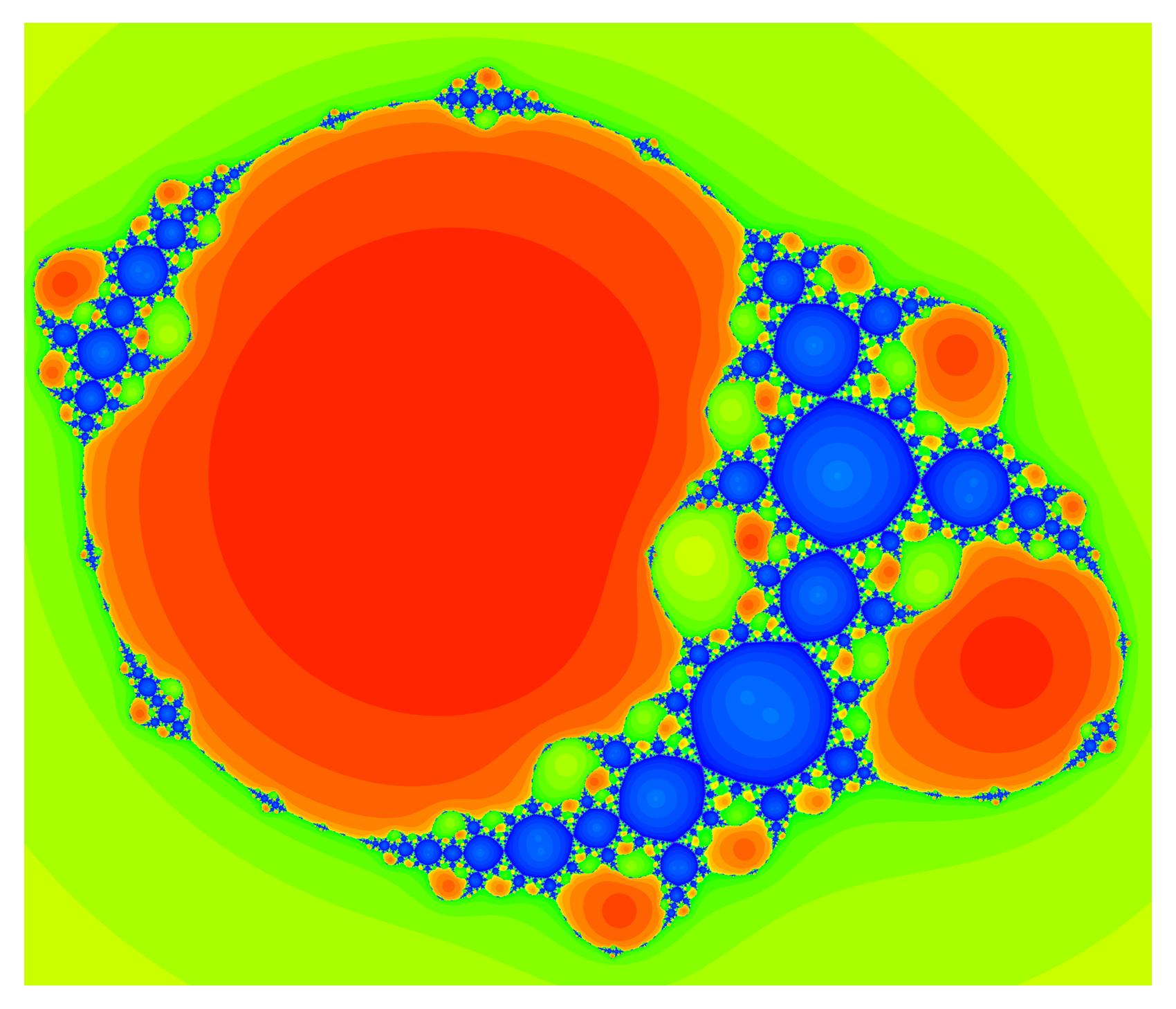}  
\includegraphics[width=0.44\textwidth, height=0.2\textheight]{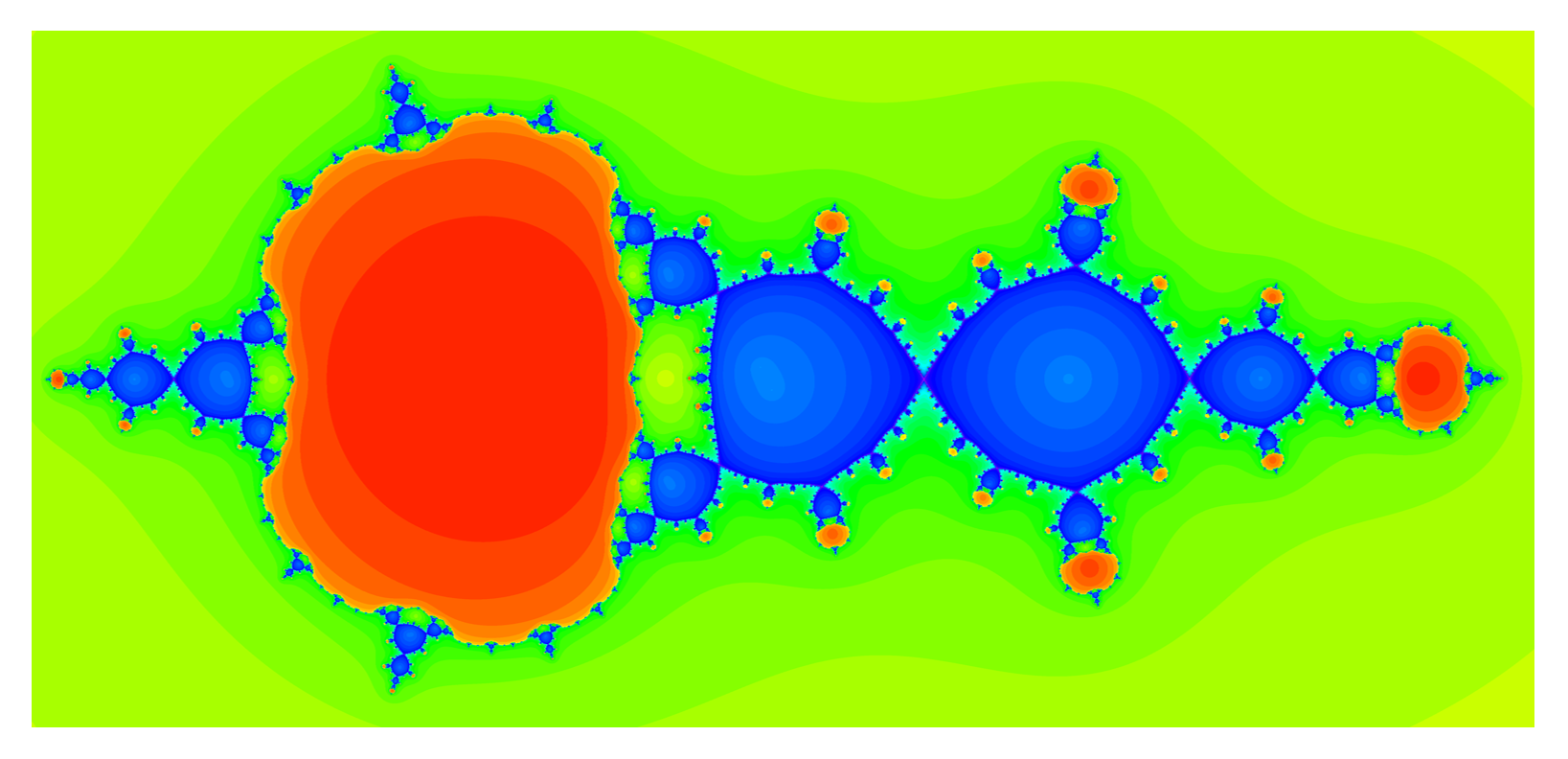}  
\includegraphics[width=0.44\textwidth, height=0.2\textheight]{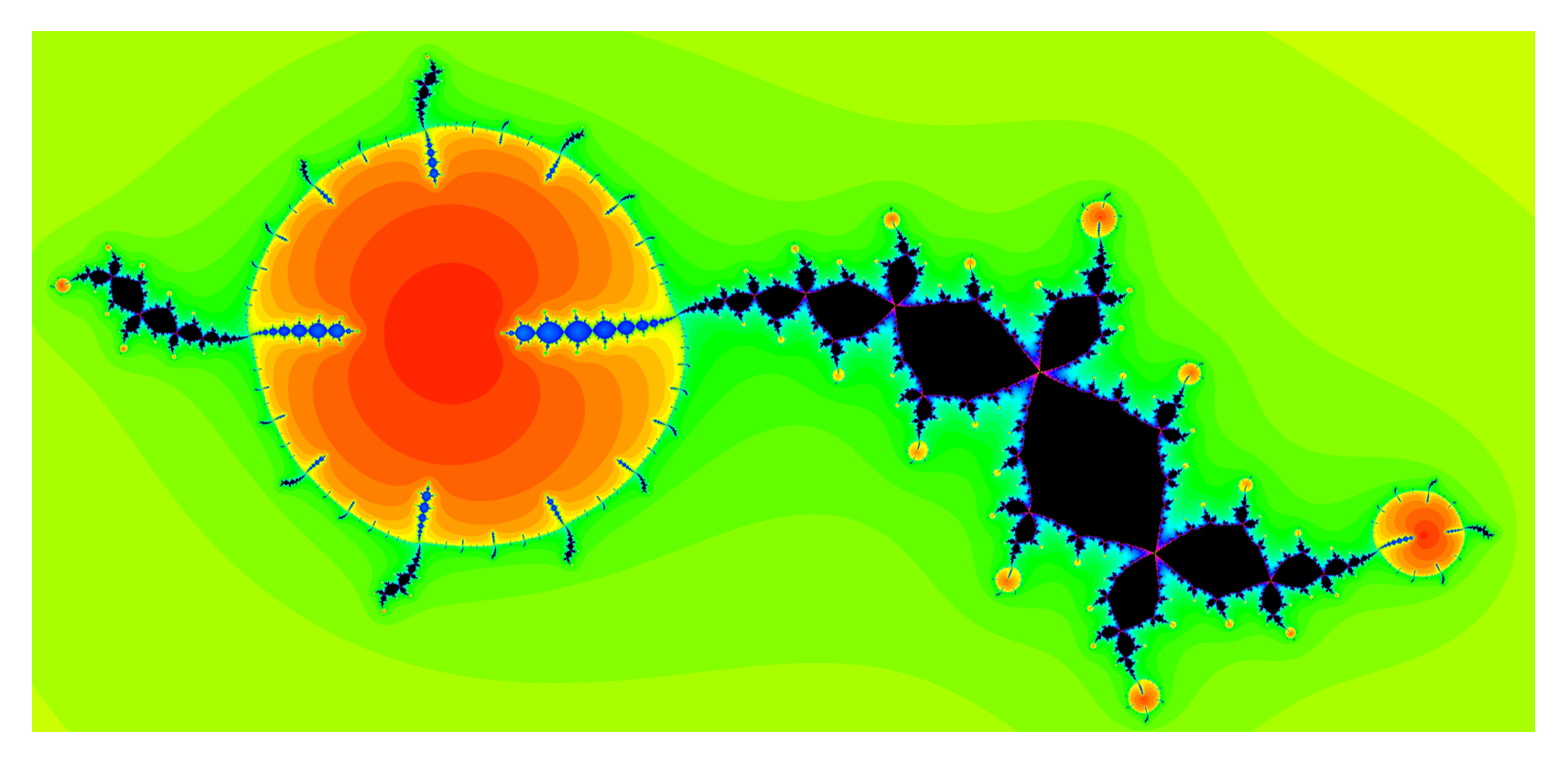}
\includegraphics[width=0.35\textwidth, height=0.2\textheight]{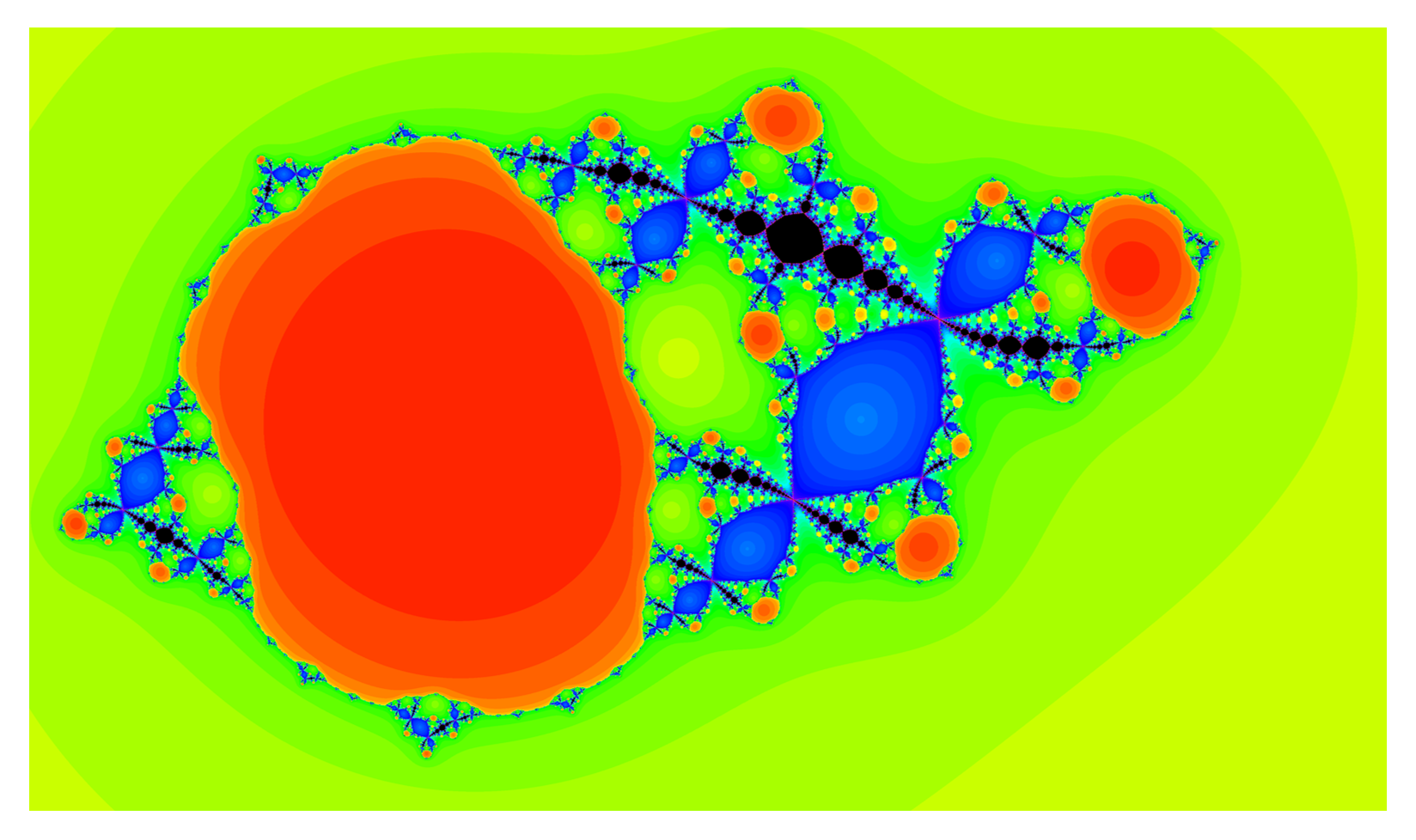} 
\caption{Representative Julia sets for parameters from components of types $A$, $B$, and $D$ in the right-hand side of Figure~\ref{fig:secper20vertd3} (in the sense of Milnor \cite{milnor20099}).}
\label{fig:milnorTypes}
\end{figure}

To support this interpretation, Figure~\ref{fig:milnorTypes} shows representative Julia sets associated with parameters selected from components of each type. 
These images reveal regions of the dynamical plane consistent with the Julia sets arising in the corresponding hyperbolic components of the slice $\mathrm{Per}_2(0)$ in the cubic polynomial family.

\begin{figure}[ht!]
   \centering   
  \includegraphics[width=0.35\textwidth, height=0.25\textheight]{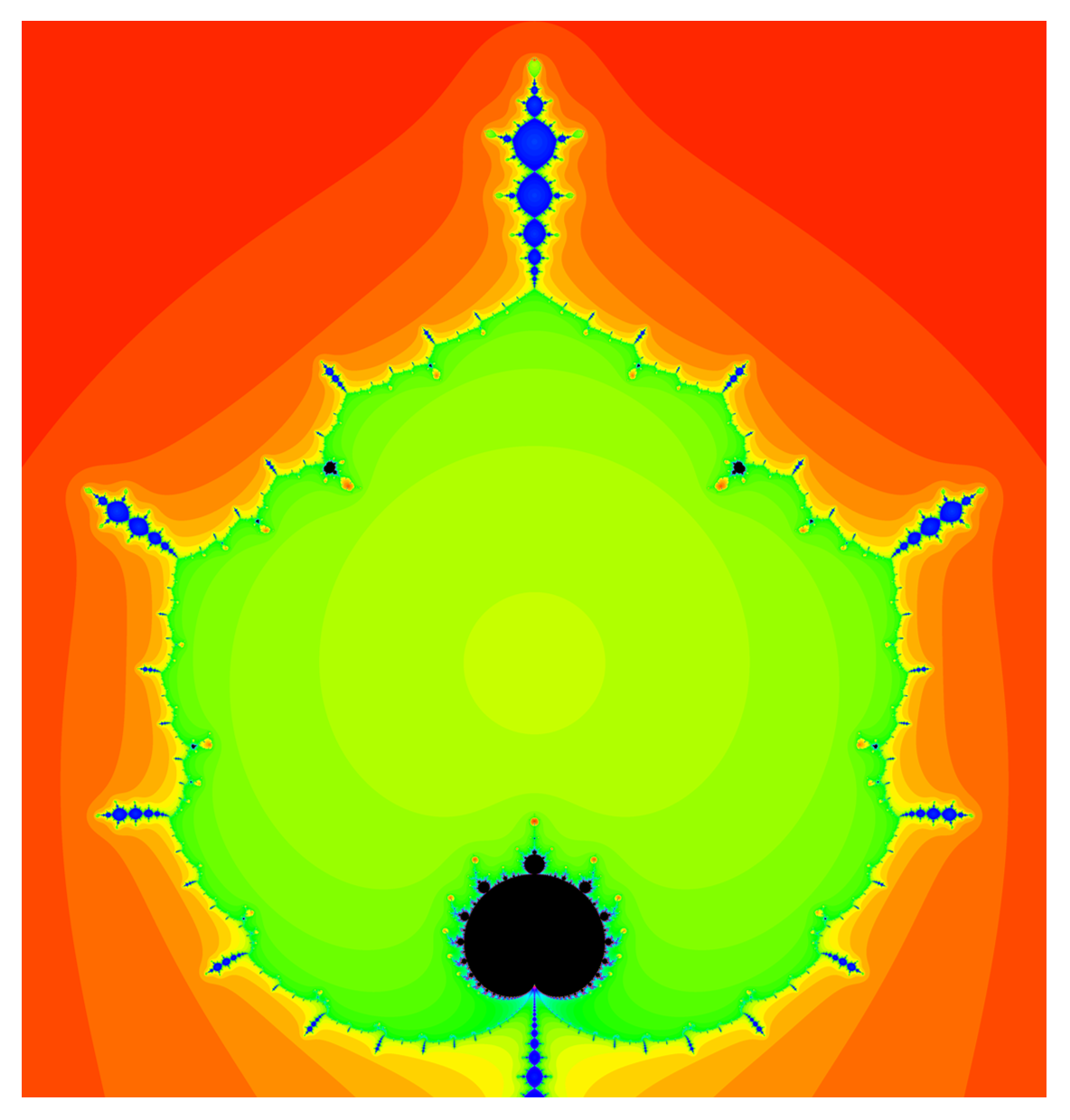} 
\includegraphics[width=0.4\textwidth, height=0.25\textheight]{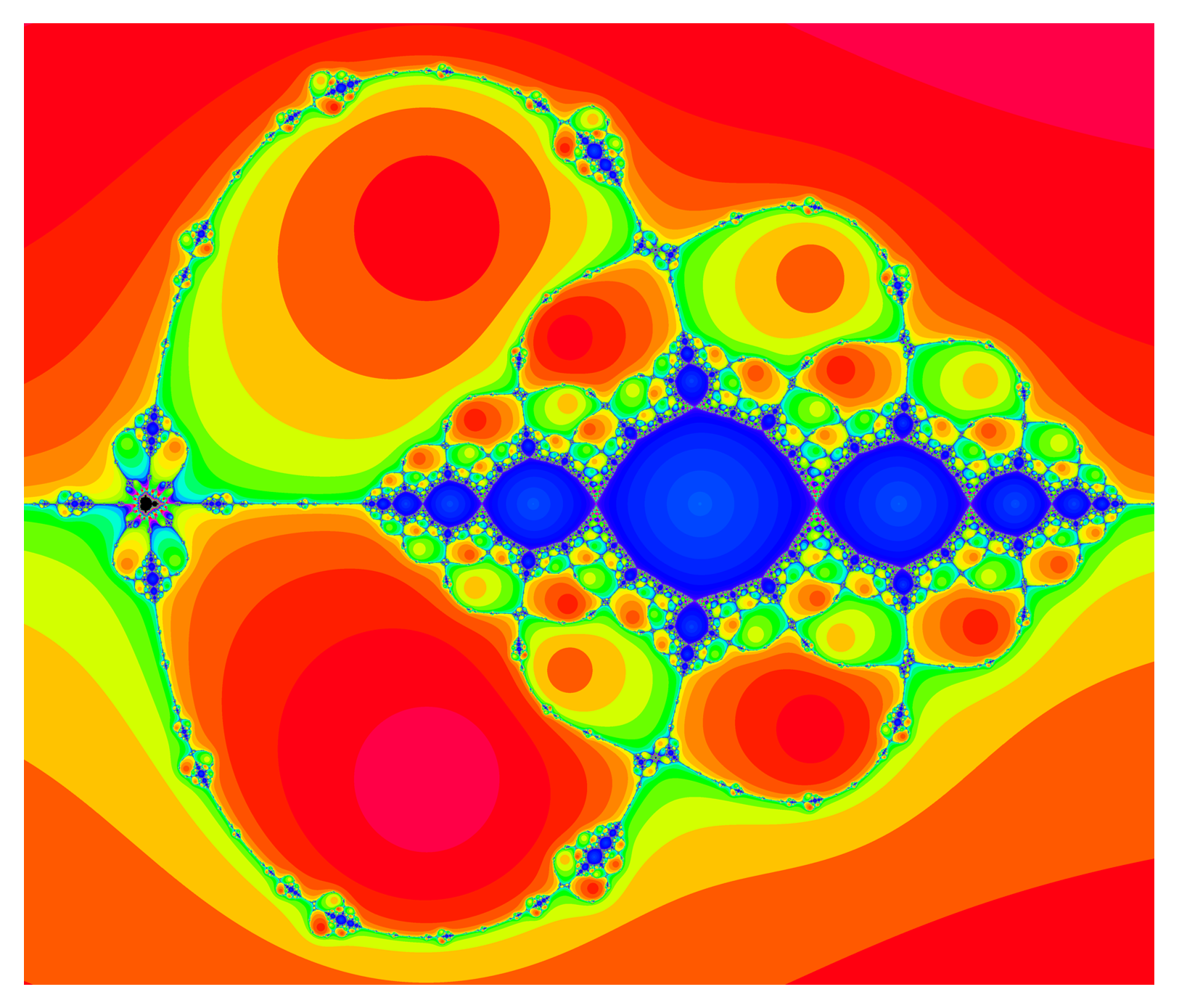} 
\caption{Zoomed views of $\mathscr{S}_2(0)$ from Figure~\ref{fig:secper20vertd3}.}
   \label{fig:per2(0)ZOOM}
\end{figure}

Figure~\ref{fig:per2(0)ZOOM} shows a zoom of Figure~\ref{fig:secper20vertd3}, where regions exhibiting structures resembling Julia sets appear embedded in the parameter space of $\mathscr{S}_2(0)$. On the left, one observes a hyperbolic component in which the free critical point is attracted to infinity after finitely many iterations. This structure has a leaf-like shape, with infinitely many Julia-like sets attached at cusps and pointing outward, as well as infinitely many baby Mandelbrot sets attached at cusps and pointing inward toward the leaf. On the right, one observes a Julia-like structure coexisting in close proximity with a Mandelbrot-like set, forming an intricate configuration, with both structures symmetric with respect to the real axis.

\subsection{Slices $\mathscr{S}_1(\lambda)$}

Fix $d$. For visualization purposes, we focus on small values of $d$, since the number of hyperbolic components increases significantly as $d$ grows. 
We investigate the slices $\mathscr{S}_1(\lambda)$ for different values of $\lambda$. Using the parametrizations derived in Subsection~\ref{subsec:lambda1}, we explore the geometry of the corresponding connectedness loci and analyze how their structure varies with $\lambda$ and $d$. 
Particular attention is given to the emergence of recurrent patterns and their relation to known features in parabolic slices of cubic polynomial parameter spaces.

We distinguish between two cases, according to the arithmetic nature of the rotation angle $\theta$: the rational case, and the more delicate case of irrational values.
\subsubsection{For $\theta \in \mathbb{Q}$}
Fixing $\lambda = e^{2\pi i \theta}$ for selected values of $\theta$, we investigate slices of the form $\mathscr{S}_1(\lambda)$ using the parametrization $f_{\beta,\lambda}$. By systematically sampling the parameter plane, we examine how the geometry depends on the rotation number and look for characteristic patterns, including parabolic bifurcations and Mandelbrot-like structures, in comparison with cubic polynomial families.

\begin{figure}[ht!]
   \centering
       \includegraphics[width=0.31\textwidth, height=0.3\textheight]{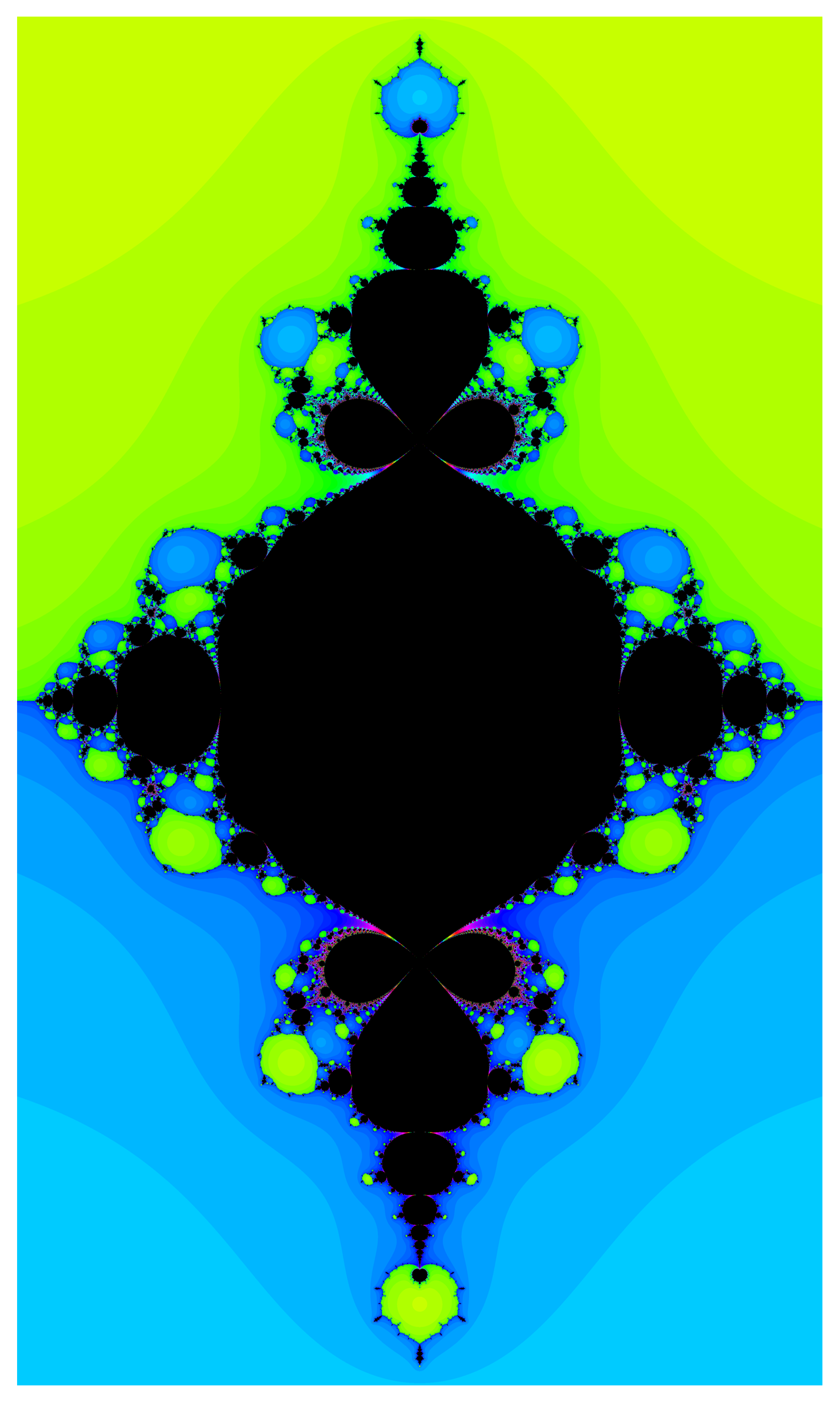} 
    \includegraphics[width=0.45\textwidth, height=0.3\textheight]{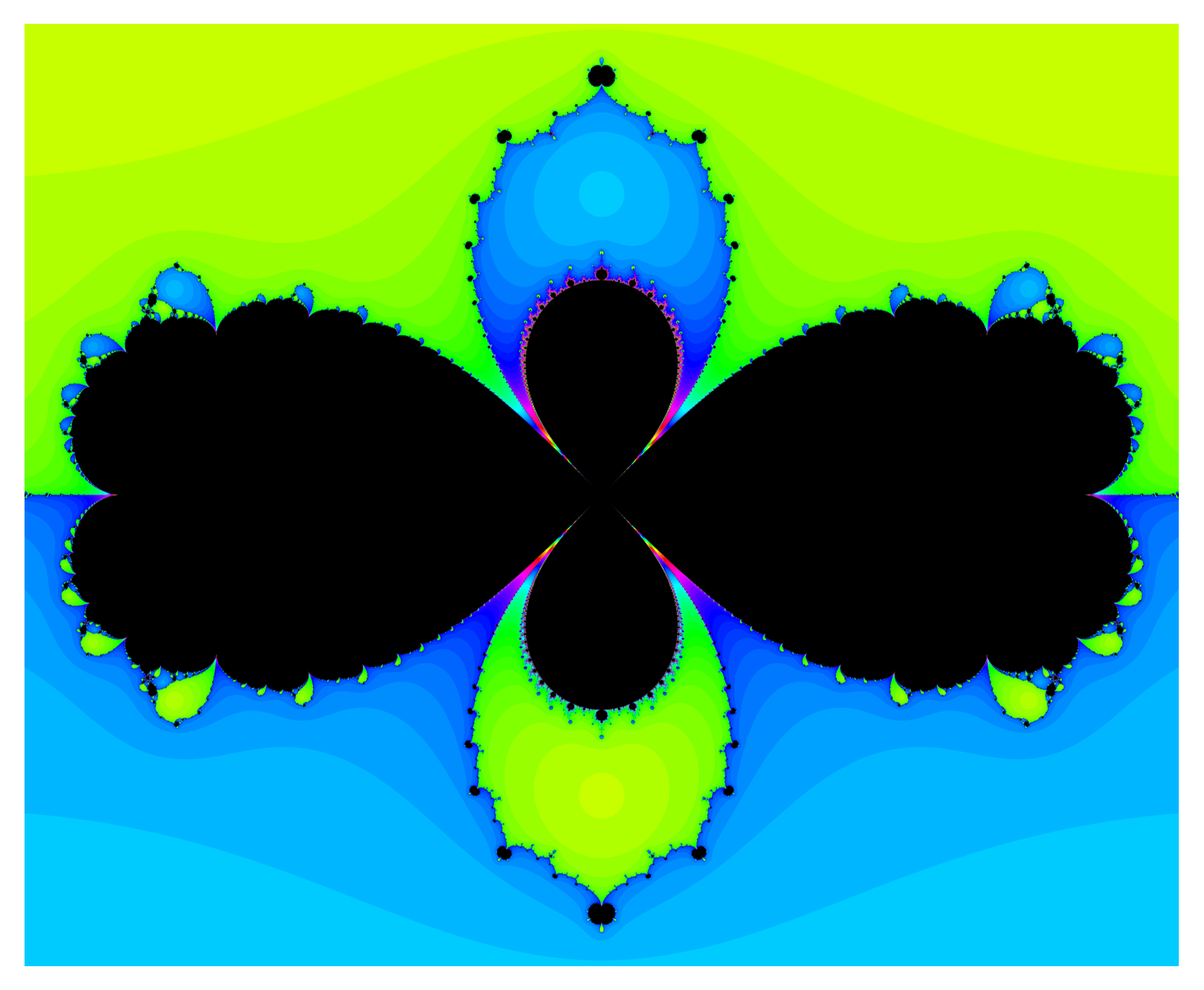} 
    \includegraphics[width=0.44\textwidth, height=0.33\textheight]{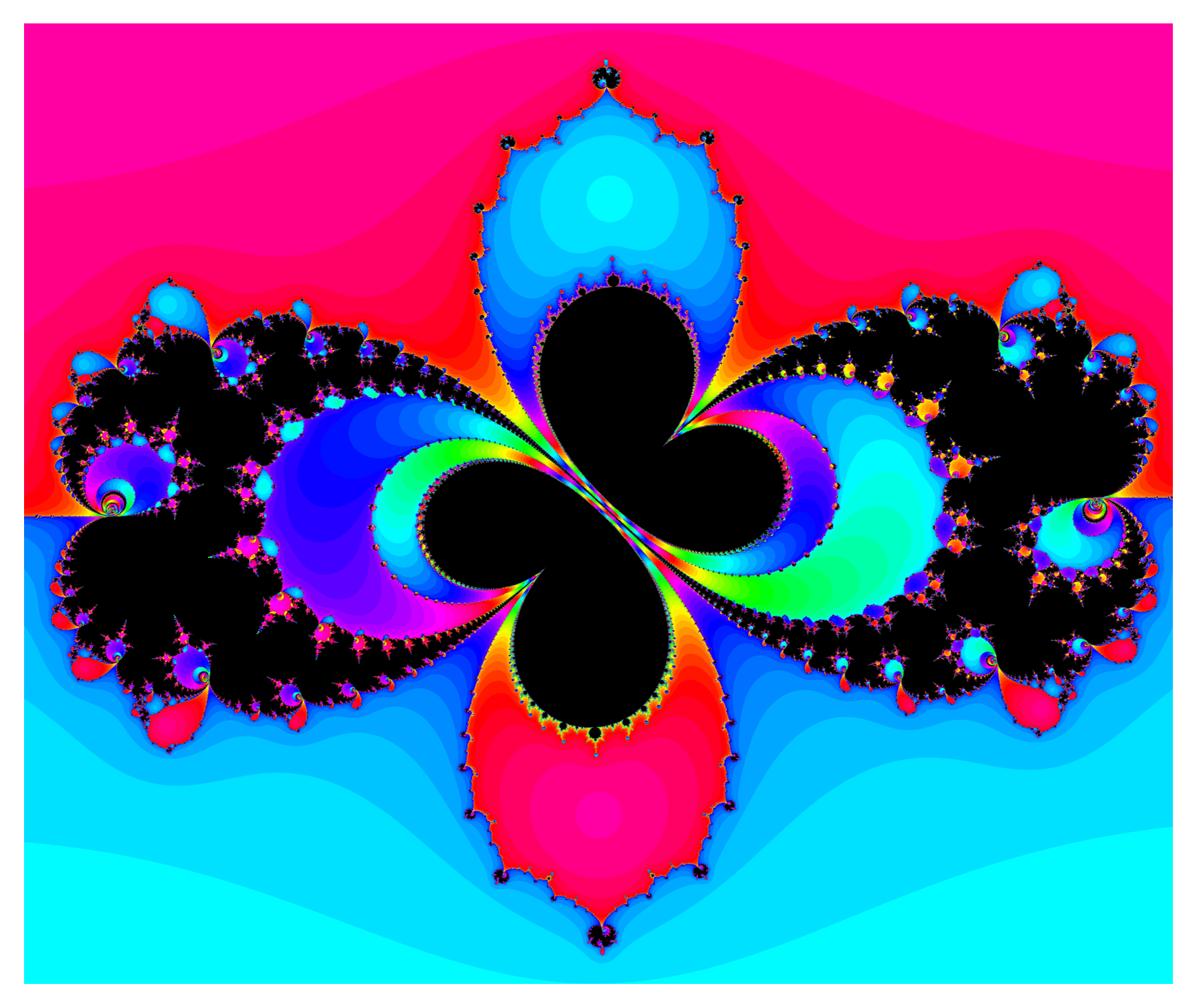}%
    \includegraphics[width=0.45\textwidth, height=0.33\textheight]{fig/plt_per1_cubic_parab_Large.png} 
\caption{Connectedness loci $\mathscr{M}_{\lambda}$ in $\mathscr{S}_1(\lambda)$ for $\lambda = e^{2\pi i \theta}$, where $\theta \in \mathbb{Q}$.} 
   \label{fig:secper1lambdauntercio}
 \end{figure}

Figure~\ref{fig:secper1lambdauntercio} suggests that the parameter slice $per_1(\lambda)$ for $\lambda=-1,1$ and $\theta=\frac{1}{180}$, in the case $d=3$, exhibits a rich parabolic bifurcation structure. The latter case can be compared with $\mathrm{Per}_1(\lambda)$ for parabolic cubic polynomials.

In particular, the observed Julia-like geometries are consistent with the local dynamics near parabolic fixed points and appear to reflect the influence of the quadratic normal form $e^{2\pi i \theta} z + z^2$. This supports the interpretation that, at small scales, the parameter space inherits universal features from lower-degree families.

A similar organization is observed in parabolic slices of cubic polynomials of the form $e^{2 \pi i \theta} z + b z^2 + z^3$, suggesting that these structures persist across degrees and may be governed by common underlying mechanisms. 
See also~\cite{zhang2025parabolic} for related phenomena.

\subsubsection{For $\theta \notin \mathbb{Q}$}

\begin{figure}[ht!]
   \centering   
   \includegraphics[width=0.45\textwidth, height=0.5\textheight]{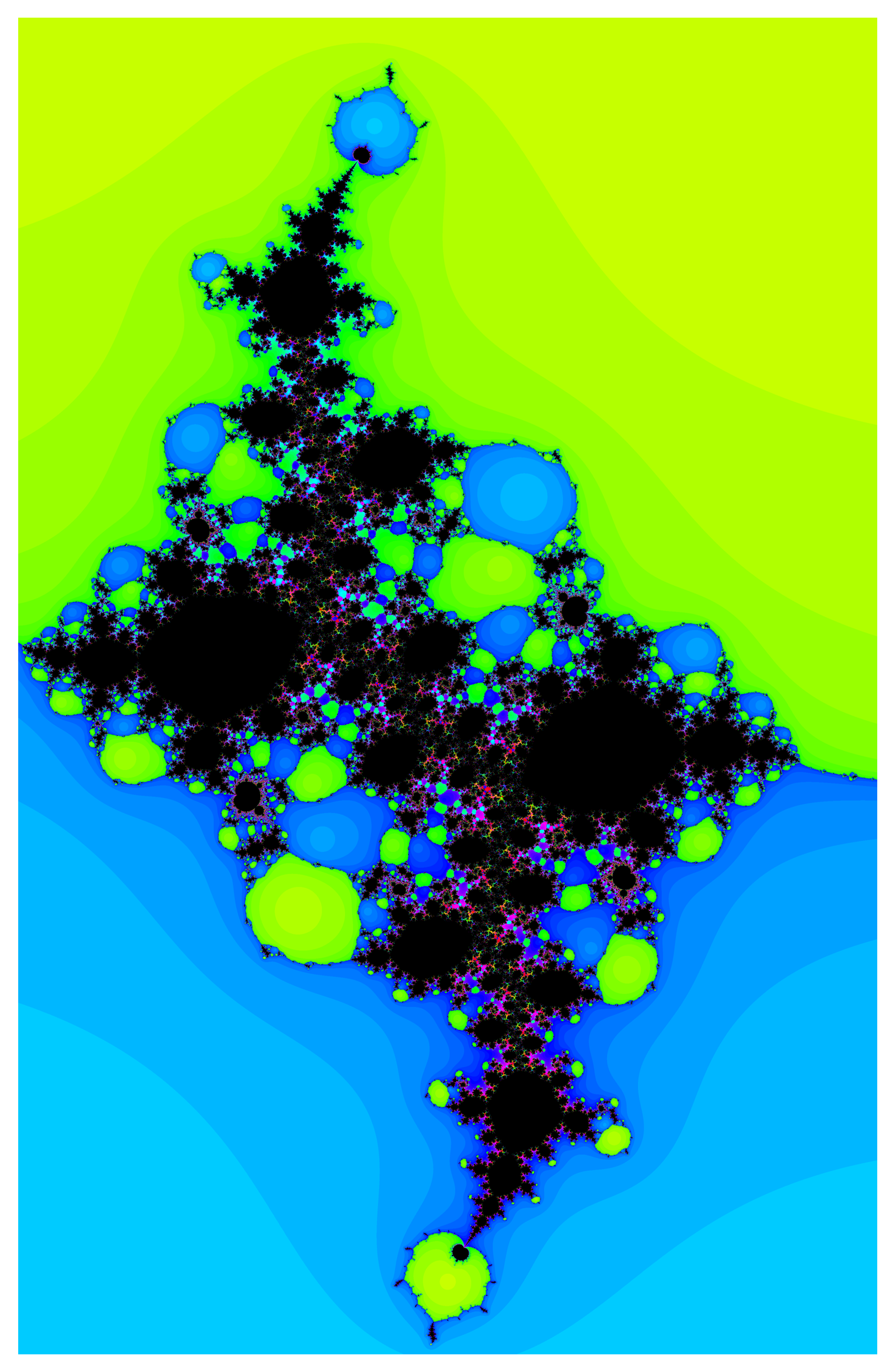}    
   \includegraphics[width=0.53\textwidth, height=0.5\textheight]{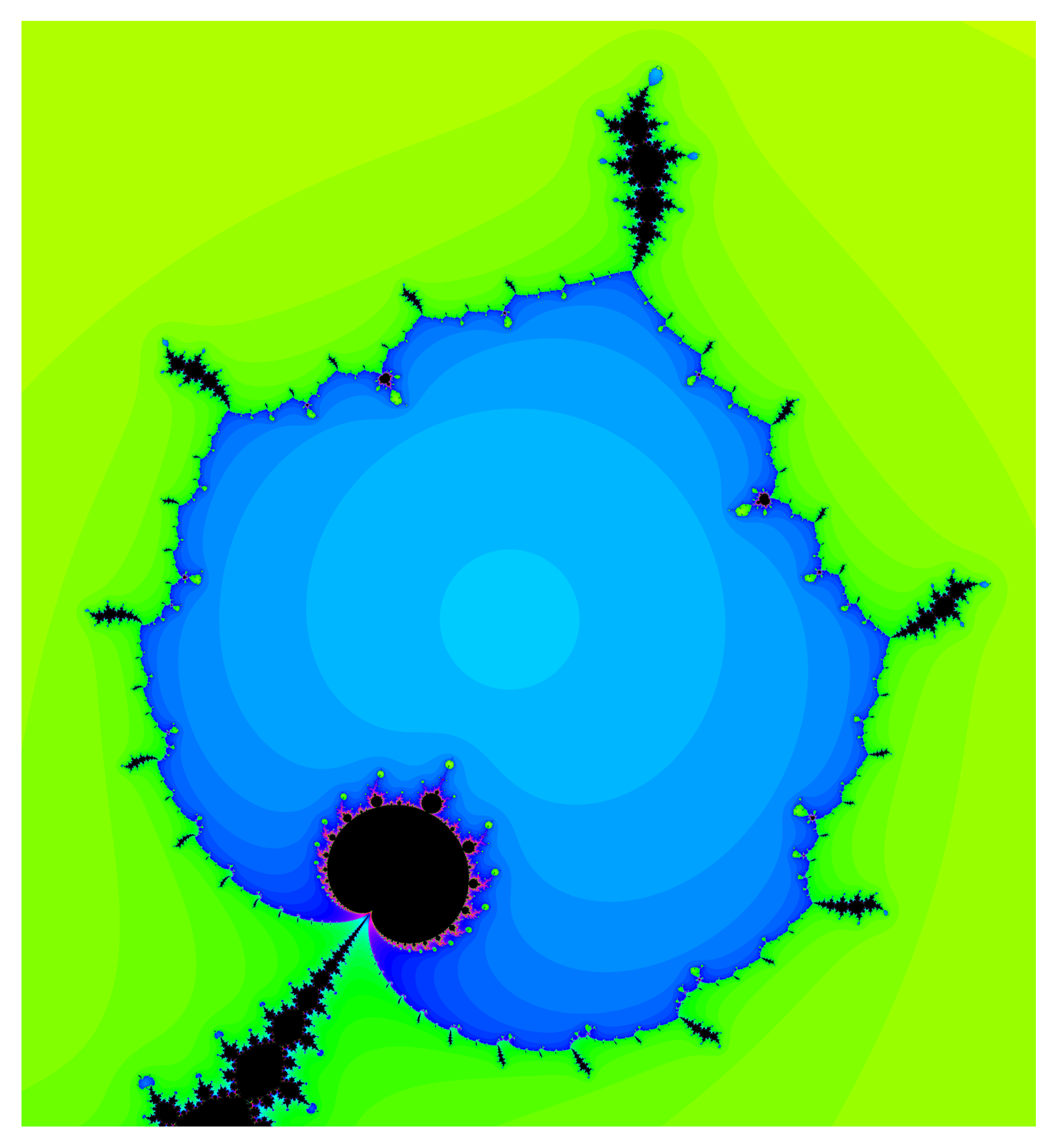} 
   \caption{Connectedness locus of $\mathscr{S}_1(\lambda)$ in the $a$-parameter plane for $d=3$ with $\theta=\frac{\sqrt{5}-1}{2}$. The parameter space exhibits intricate filamentary structures and fine-scale patterns, consistent with the complex behavior associated with irrational rotation numbers satisfying the Brjuno condition.}
   \label{fig:golden}
\end{figure}

We consider the inverse of the golden ratio, 
$ \theta=\frac{\sqrt{5}-1}{2}.$
Figure~\ref{fig:golden} shows the connectedness locus of $\mathscr{S}_1(\lambda)$. 

This choice is motivated by the fact that $\theta$ is a classical example in holomorphic dynamics and a distinguished Brjuno number, for which linearization phenomena are expected to occur near neutral fixed points. In particular, it is associated with the well-known Douady conjecture \cite{douady1987disques}, which remains open in the setting of rational maps.

From a dynamical perspective, rotation numbers satisfying the Brjuno condition are associated with linearizable behavior near neutral cycles. In parameter space, this is often reflected in the presence of finely structured regions. The numerical experiments presented here reveal intricate geometric patterns that are consistent with this general phenomenon.

In contrast to the rational case, the geometry of the parameter space exhibits a marked loss of symmetry and a significantly higher level of complexity. In particular, one observes filamentary structures that appear to organize along curves separating regions with distinct dynamical behavior, especially between bounded and escaping critical orbits. This suggests that the global geometry is strongly influenced by a delicate interaction between critical orbit relations and the arithmetic properties of $\theta$.

At finer scales, the parameter space displays regions with a highly textured appearance, reminiscent of Julia sets. This indicates that the boundary of the connectedness locus may encode dynamical features typically associated with invariant sets in the dynamical plane. While a precise identification of these structures remains beyond the scope of the present work, their presence is consistent with phenomena observed in polynomial families with neutral dynamics.

Accessing conjectural objects such as the Jordan curves in the bifurcation loci arising in the work of Zakeri~\cite{zakeri1999dynamics} through numerical methods appears to be a challenging problem. To the best of our knowledge, no effective algorithm is currently available for visualizing such objects at high resolution in this setting, much less in more general contexts such as the one considered in this work. For this reason, we restrict our attention to qualitative features that can be reliably observed at the computational scale considered here.

In analogy with parameter spaces of cubic polynomials with a persistent neutral or parabolic point, one expects the presence of Mandelbrot-like structures organized by critical orbit relations. Nevertheless, these structures may be highly nontrivial: in some families, Mandelbrot-like sets are known to be non-computable~\cite{coronel2018non}. This suggests that the geometry observed numerically may only reflect part of a richer underlying structure.

\subsection{Slices $\mathscr{S}_2(\lambda)$}
Fixing $\lambda = e^{2\pi i \theta}$ and varying both $\theta$ and the degree $d$, we carry out numerical experiments to investigate the geometry of the connectedness locus $\mathcal{M}(\theta)$.

Our goal is to compare the observed structures with the period-two slices $\mathrm{Per}_2(\lambda)$ in the parabolic setting of cubic polynomials, in order to provide evidence supporting the conjectural framework described in the introduction.

\subsubsection{For $\theta \in \mathbb{Q}$}
\FloatBarrier

\begin{figure}[ht!]
   \centering   
   \includegraphics[width=0.45\textwidth, height=0.3\textheight]{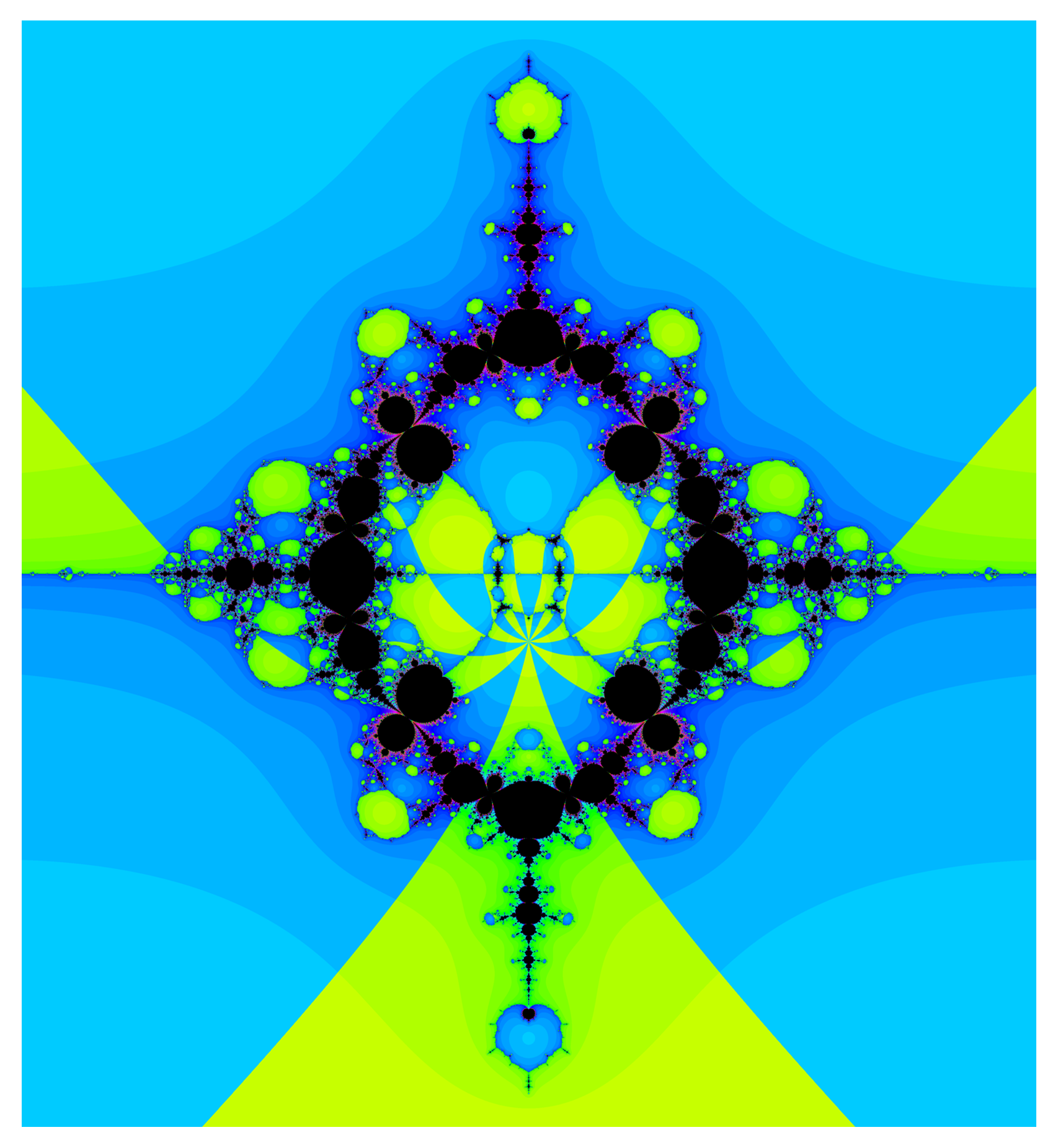} 
   \includegraphics[width=0.45\textwidth, height=0.3\textheight]{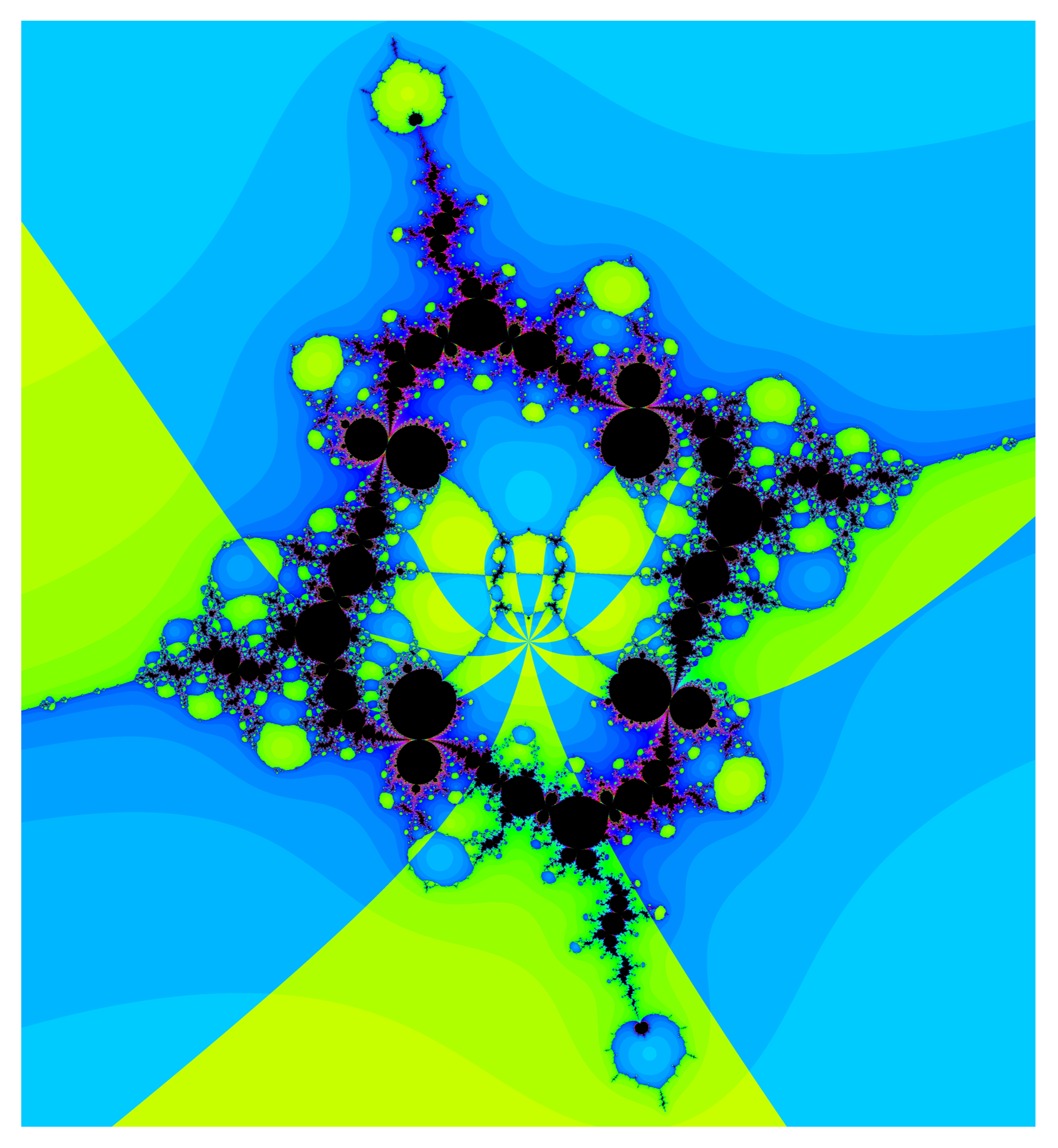} 
   \includegraphics[width=0.45\textwidth, height=0.3\textheight]{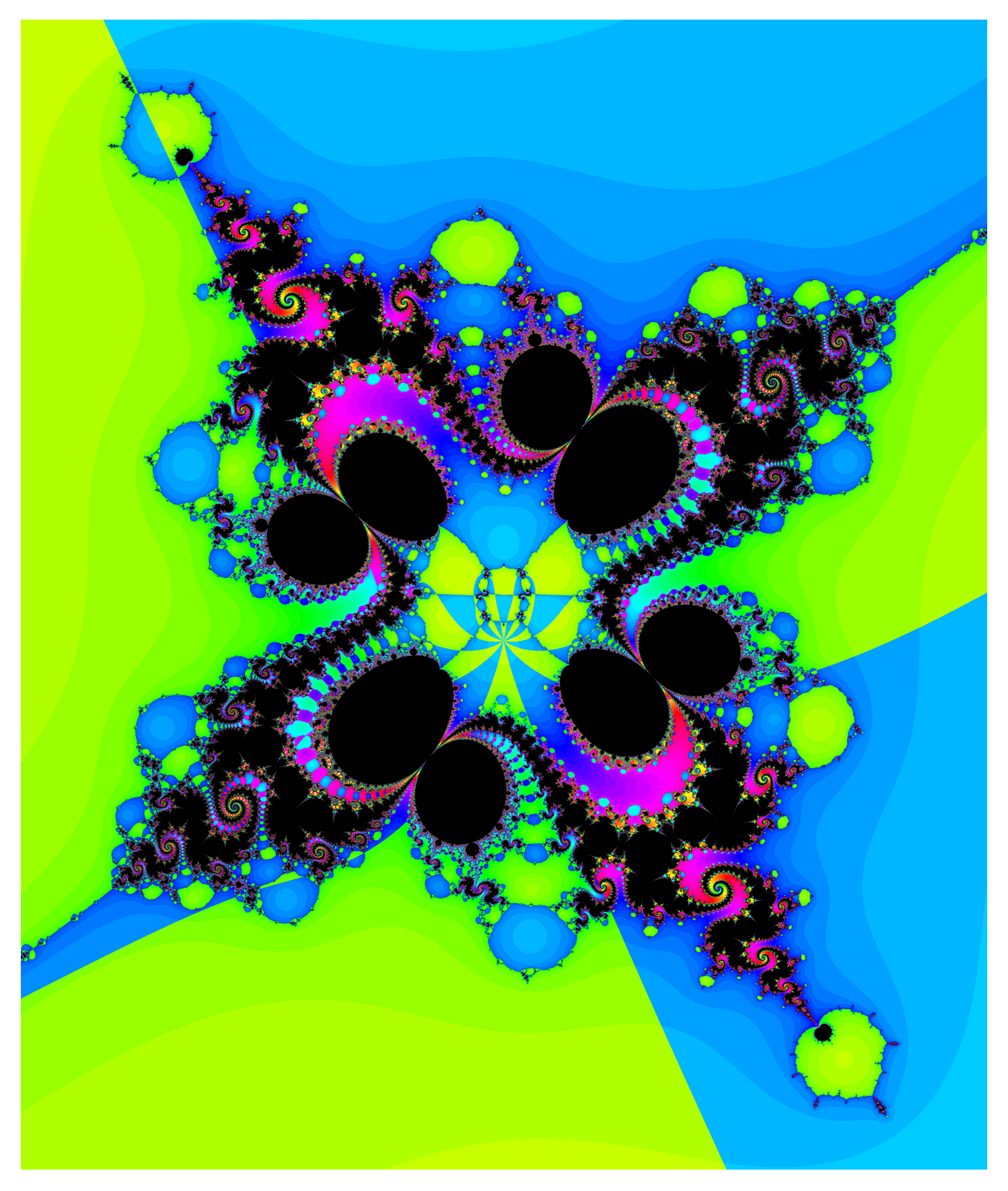} 
   \includegraphics[width=0.45\textwidth, height=0.3\textheight]{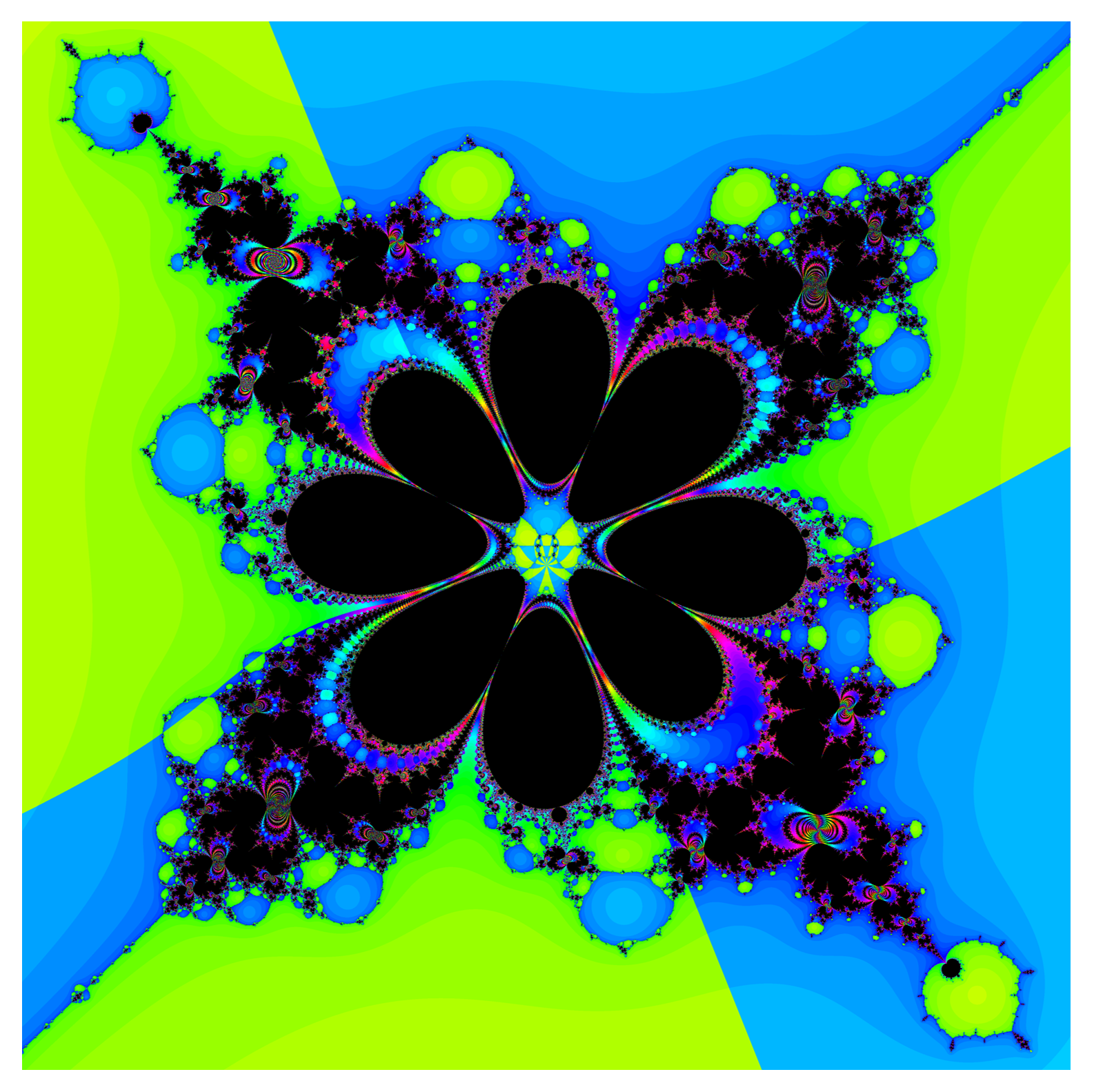}  
   \caption{Connectedness loci $\mathcal{M}_{\theta}$ for $\theta = \frac{1}{2}, \frac{1}{3}, \frac{1}{11}$ and $\frac{1}{79}$. The parameter space exhibits a decomposition into structured regions consistent with hyperbolic components, together with increasingly intricate boundaries as the denominator of $\theta$ grows.}
   \label{fig:GalleryS2}
\end{figure}

Figure~\ref{fig:GalleryS2} shows the connectedness loci $\mathcal{M}_{\theta}$ for increasing values of $q$, starting from small values and progressing to larger ones. Across these parameter slices, we observe persistent geometric structures that closely resemble the period-two slices $\mathrm{Per}_2(\lambda)$ in the parameter space of cubic polynomials.

A prominent feature of these slices is the organization of the parameter space into well-defined regions that appear to correspond to hyperbolic components. These regions are arranged around distinguished parameters, which act as organizing centers for the surrounding geometry. This behavior is consistent with the presence of attracting cycles capturing one or both critical orbits, and reflects the fundamental role of critical orbit relations in structuring the parameter space.

As $\theta = p/q$ with $q \to \infty$, the parameter space becomes increasingly intricate. In particular, the boundary of the connectedness locus develops thin filamentary structures and small-scale features that accumulate near central regions. This suggests an enrichment of the parabolic bifurcation structure, analogous to phenomena observed in cubic polynomial families; compare with \cite{zhang2025parabolic}.

At finer scales, one observes localized structures reminiscent of Mandelbrot-like sets, which appear to be organized by repeated critical orbit configurations. This indicates that mechanisms similar to polynomial-like renormalization may be present in this rational setting, even though a precise formulation remains open.

\begin{figure}[ht!]
   \centering   
   \includegraphics[width=0.3\textwidth, height=0.3\textheight]{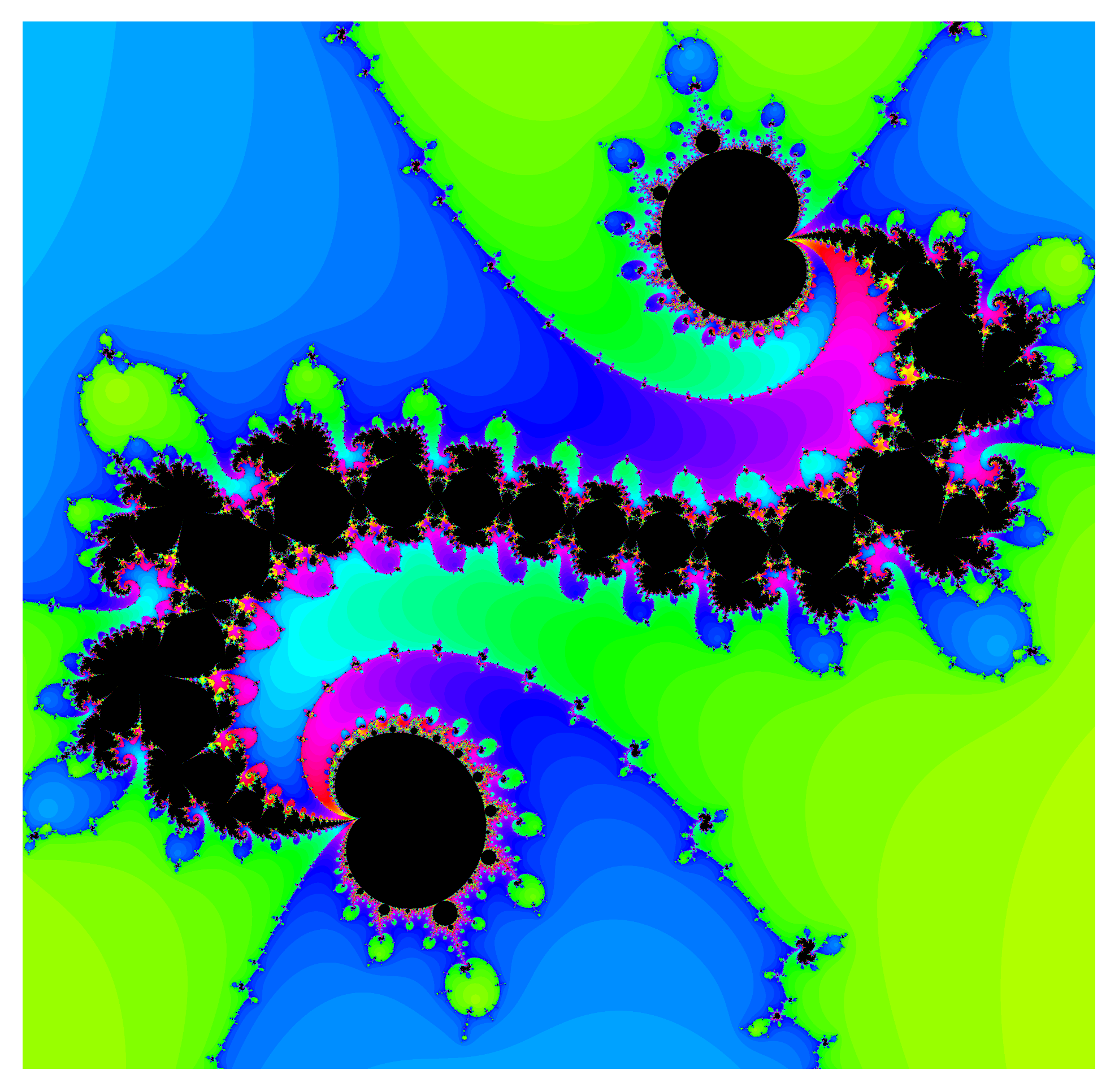} 
   \includegraphics[width=0.3\textwidth, height=0.3\textheight]{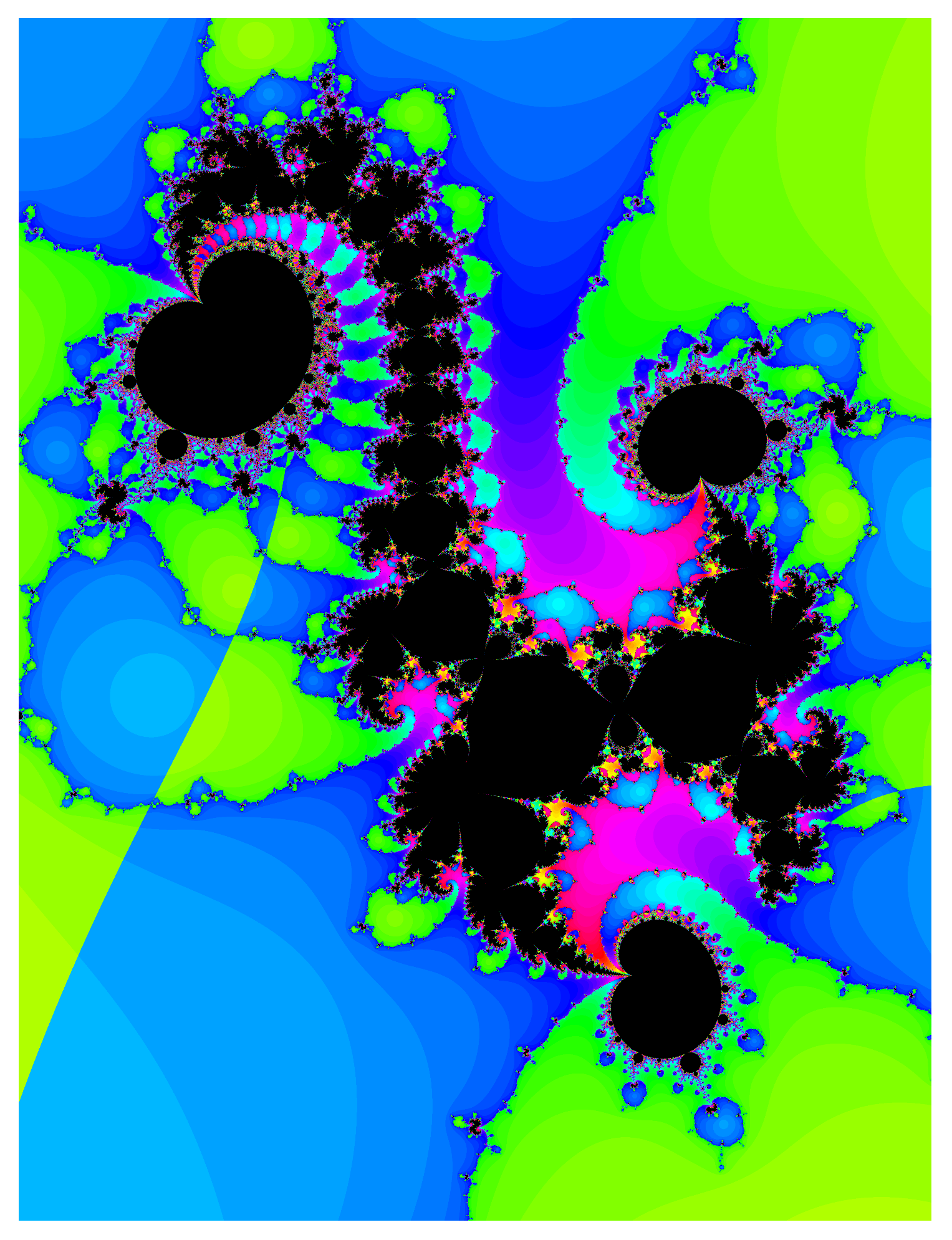} 
   \includegraphics[width=0.35\textwidth, height=0.3\textheight]{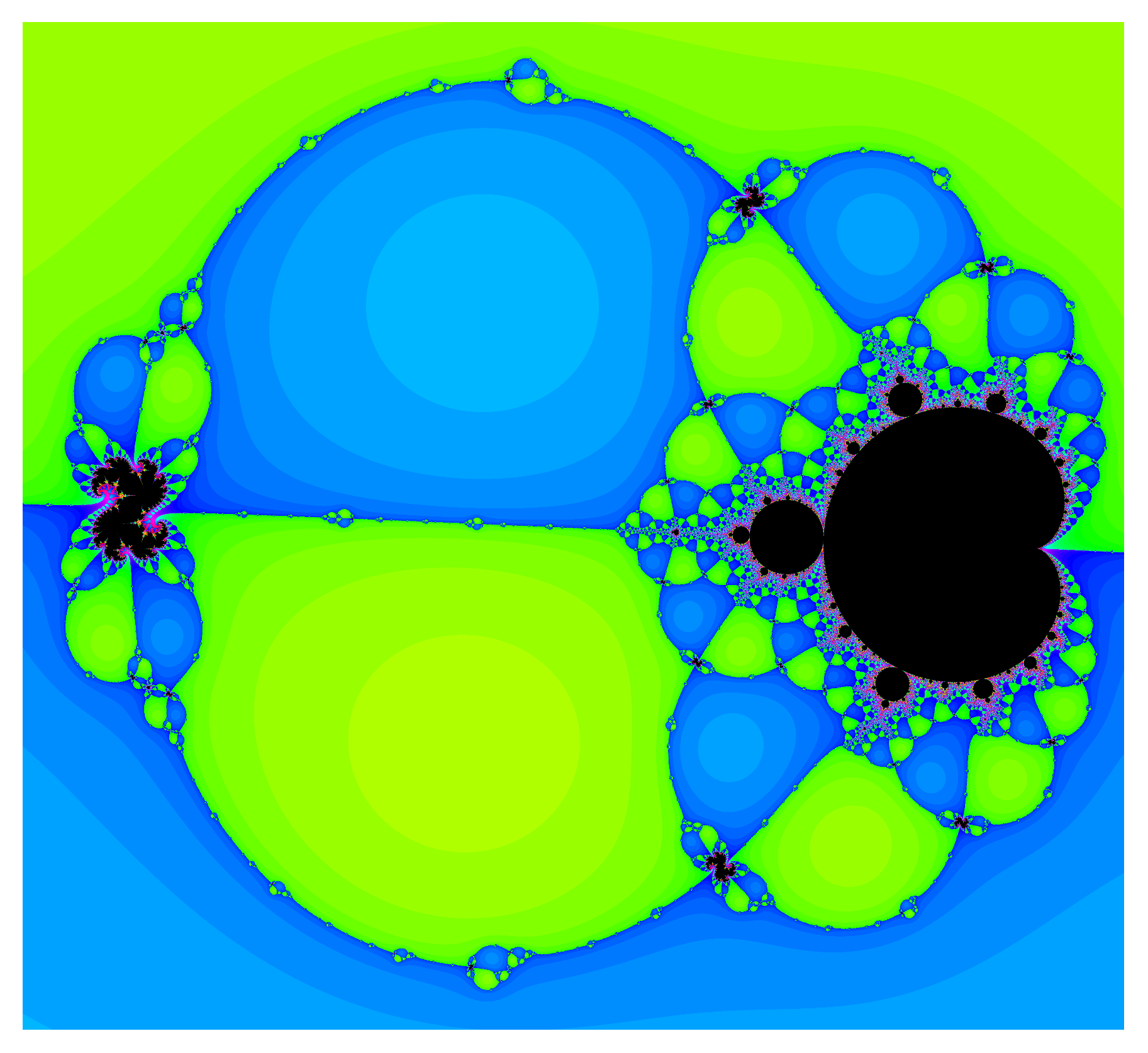} 
 \caption{Zoomed-in views illustrating the fine-scale structure of the parameter space in Figure~\ref{fig:GalleryS2}.}
   \label{fig:per2_d111ZOOM}
\end{figure}

Figure~\ref{fig:per2_d111ZOOM} shows zoomed-in views of regions near the center of the parameter plane in the rightmost panel of Figure~\ref{fig:GalleryS2}. These reveal a hierarchical organization of components, with smaller structures attached to larger ones in a manner consistent with iterated bifurcation processes. In particular, one observes the coexistence of distinct dynamical regimes corresponding to different configurations of the critical orbits, including regions resembling $\mathrm{Per}_1(\lambda)$ for parabolic cubic polynomials. In the rightmost panel, one also observes a Mandelbrot set coexisting with a parabolic Julia set associated with the family $\lambda z + z^2$.

For an intuitive comparison, see Figure~\ref{fig:secper20vertd3}: the blue component may be viewed as being replaced by a parabolic Julia set occupying the same region.

\subsubsection{For $\theta \notin \mathbb{Q}$}
\FloatBarrier

\begin{figure}[ht!]
   \centering
   \includegraphics[width=0.7\textwidth, height=0.55\textheight]{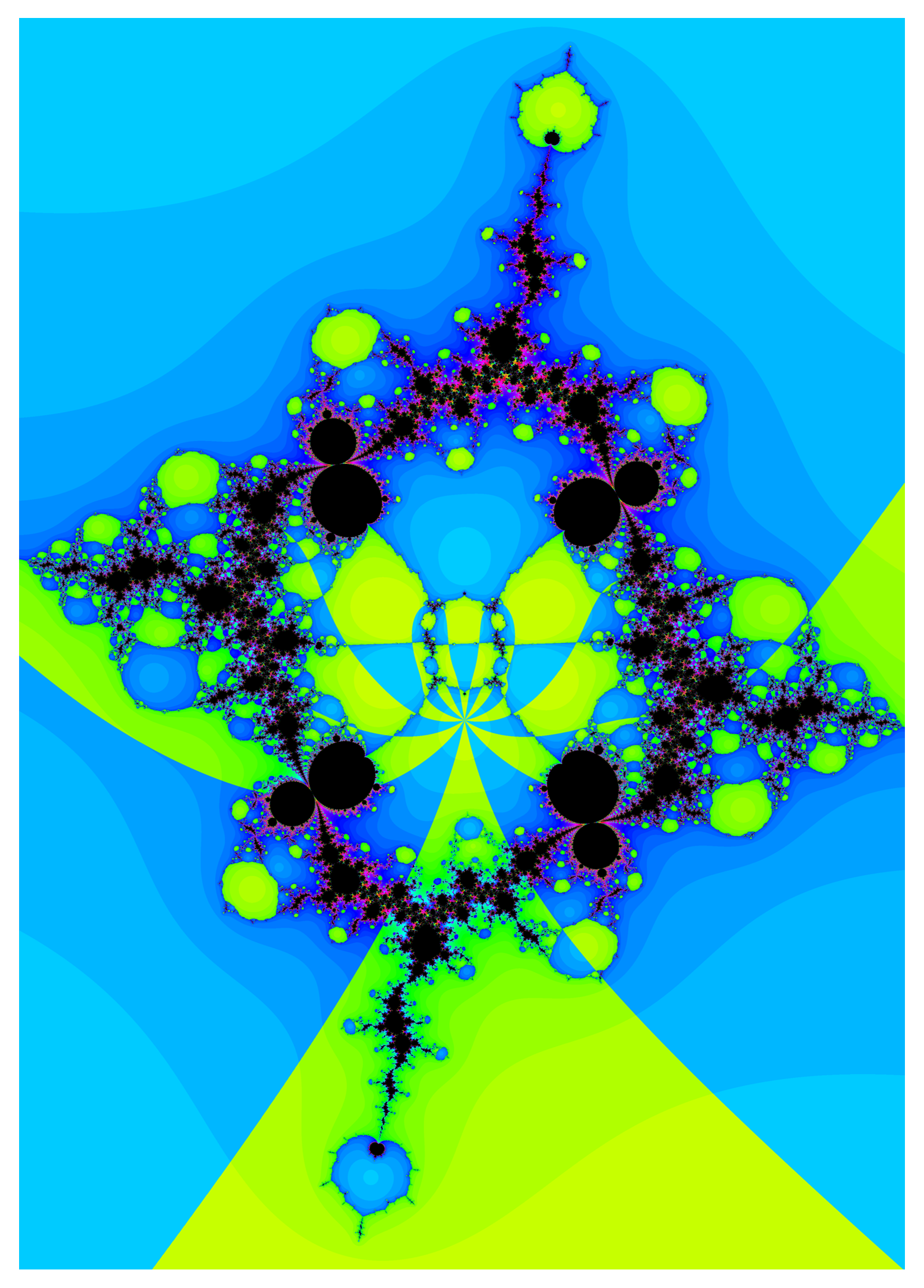} 
\caption{Connectedness locus $\mathcal{M}_{\theta}$ in the slice $\mathscr{S}_2(\theta)$ for $d=3$ and $\theta=\frac{\sqrt{5}-1}{2}$.}
   \label{fig:secper2lambda11}
\end{figure}

Figure~\ref{fig:secper2lambda11} illustrates the connectedness locus $\mathcal{M}(\theta)$ in $\mathscr{S}_2$ for the inverse golden ratio $\theta=\frac{\sqrt{5}-1}{2}$. Numerical experiments reveal intricate geometric structures, including filamentary patterns organized along curves that resemble the loci $\mathrm{Per}_1(\lambda)$ and $\mathrm{Per}_2(\lambda)$ observed in parabolic slices of cubic polynomial families (see Figure~\ref{fig:golden}). These features suggest the presence of analogous organizing mechanisms in the present setting, although their precise dynamical interpretation remains unclear. The topological structure of the slice $\mathscr{S}_{2}(\lambda)$, as well as its relationship with $\mathrm{Per}_2(\lambda)$ for cubic polynomial families, is still largely unexplored and constitutes a promising direction for future research.

\subsection{Shrimp-like structures in the Blaschke slice}
We now focus on a particular subfamily of \eqref{baranskifamily} where additional structure becomes visible, namely the Blaschke slice. This setting provides a natural framework for observing organized regions of stability, as well as more intricate patterns such as shrimp-like structures (a term introduced in other contexts; see \cite{gallas1993structure,gallas2010structure}), arising from the interaction between the critical orbits.

This picture is closely related to the behavior observed in the Blaschke family studied in~\cite{canela2015family}, which arises as a special case of~\eqref{maineq}. While the geometry of the tongue regions differs from that of classical shrimp-like structures, both exhibit domains where attracting periodic dynamics persists, organized by bifurcation loci.

For $\vert a \vert = \vert b \vert =1$, the subfamily~\eqref{baranskifamily} consists of Blaschke products $f_{a,b}$ preserving the unit circle $\mathbb{S}^1$ and inducing degree-$d$ coverings of circle maps. In this setting, both critical points lie on $\mathbb{S}^1$ and play a central role in determining the dynamics.

We parametrize this slice by $(\omega_1,\omega_2)\in\mathbb{T}^2$, where $\omega_1=\arg(a)$ and $\omega_2=\arg(b)$, and classify parameters according to the asymptotic behavior of the critical orbits. A parameter is said to belong to a resonance region if at least one critical orbit converges to an attracting periodic cycle on $\mathbb{S}^1$.

\begin{figure}[ht!]
\centering
\includegraphics[width=0.7\textwidth, height=0.4\textheight]{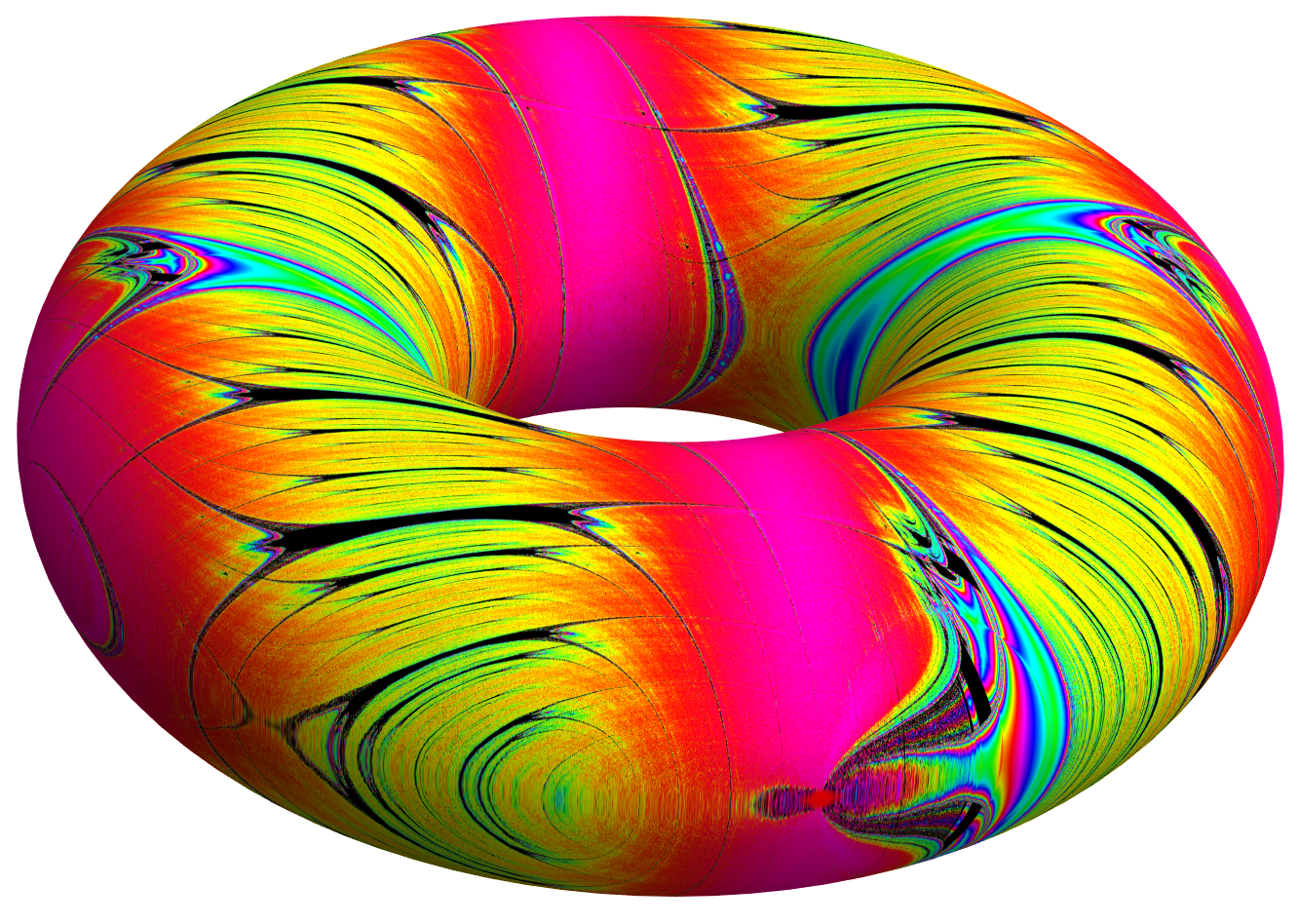} 
\caption{Parameter space of the Blaschke slice $(|a|=|b|=1)$, parametrized by $(e^{2\pi i\omega_1},e^{2\pi i\omega_2})\in\mathbb{T}^2$.}
\label{fig:Tongues2}
\end{figure}

In this parameter space, one observes tongue-like regions corresponding to attracting periodic dynamics. These tongues organize the parameter plane and reflect resonance phenomena associated with non-invertible circle coverings. In particular, they provide a natural skeleton for the global structure of the parameter space, along which more refined features develop (see Figure~\ref{fig:Tongues2}).

At finer scales, shrimp-like structures emerge within and along these tongue regions. These configurations exhibit a regular internal organization: each region corresponds to parameters for which an attracting periodic orbit persists, and is typically bounded by bifurcation curves such as saddle-node and period-doubling loci. This suggests that shrimp-like regions can be interpreted as stability windows embedded within the resonance structure defined by the tongues.

A notable feature of these structures is their hierarchical organization. Smaller shrimp-like regions appear attached to larger ones, indicating the presence of repeated bifurcation mechanisms and suggesting a form of self-similarity in parameter space. This behavior is consistent with phenomena observed in other dynamical systems, where shrimp-like sets organize stability domains within more complex regimes.

\begin{figure}[ht!]
\centering
\includegraphics[width=0.9\textwidth, height=0.3\textheight]{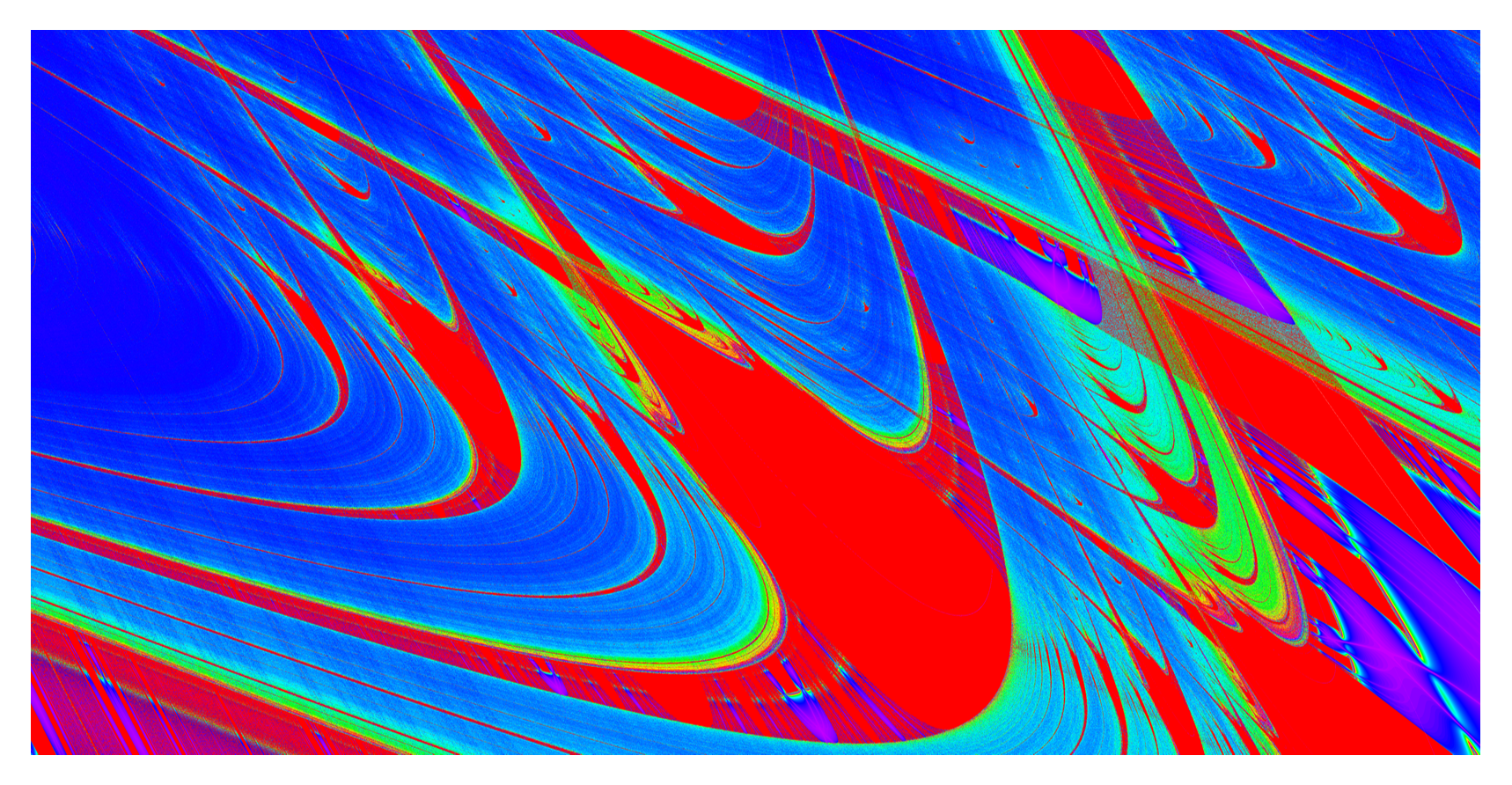}
\caption{Zoom of a shrimp-like structure in $(\omega_1,\omega_2)$-parameter plane, illustrating its internal organization and the presence of bifurcation curves delimiting stability regions.}
\label{fig:shrimp12}
\end{figure}

In the present setting, these structures arise from the interaction between the two critical orbits. This is particularly transparent in the Blaschke slice, where both critical points lie on $\mathbb{S}^1$ and directly influence the formation of resonance regions and their internal bifurcation structure. An example of a shrimp-like configuration is shown in Figure~\ref{fig:shrimp12}.

\begin{figure}[ht!]
\centering
\includegraphics[width=0.65\textwidth, height=0.4\textheight]{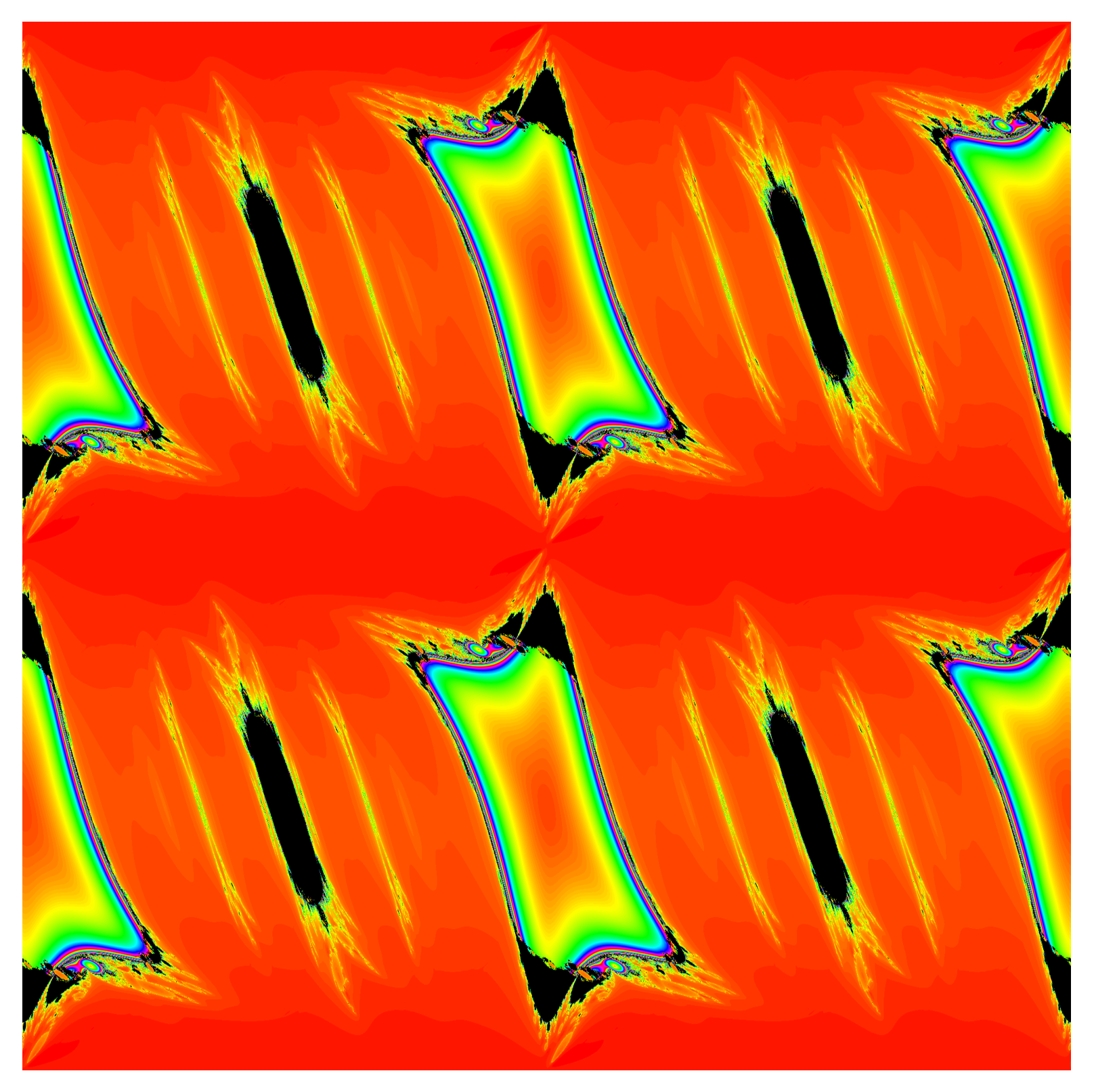}
\caption{Periodic island for $|a|\approx 1$ and $|b|\approx 1$, illustrating an isolated stability region disconnected from the main structures.}
\label{fig:island}
\end{figure}

Finally, in the region $|a|<1$ and $|b|<1$ (more precisely, writing $a=|a|e^{2\pi i \omega_1}$ and $b=|b|e^{2\pi i \omega_2}$), we observe periodic islands that are disconnected from the main Mandelbrot-like structures (Figure~\ref{fig:island}) in the $(\omega_1,\omega_2)$-parameter space. These islands correspond to isolated stability domains and further illustrate the richness of the parameter space beyond the primary resonance regions.

\section{Concluding remarks}

The numerical experiments presented in this paper reveal several recurring features across the parameter slices associated with the family~\eqref{maineq}. In both $\mathscr{S}_1(\lambda)$ and $\mathscr{S}_2(\lambda)$, we observe the emergence of Mandelbrot-like structures near singular parameters, as well as regions whose geometry closely resembles that of cubic polynomial parameter spaces. These similarities suggest that parts of the rational parameter slices may exhibit polynomial-like behavior, although a precise formulation of this correspondence remains open.

A central theme throughout this work is the dependence of the geometry on the rotation number $\theta$. Rational values tend to produce more symmetric and organized structures, often consistent with hyperbolic components and well-defined bifurcation patterns. In contrast, irrational values lead to significantly more intricate geometries, with filamentary boundaries and fine-scale structures reflecting the sensitivity of the dynamics to arithmetic properties, in analogy with classical phenomena in holomorphic dynamics.

The Blaschke slice provides a complementary perspective, where the parameter space is organized by resonance tongues associated with attracting periodic orbits on the unit circle. Within this framework, shrimp-like structures emerge as stability windows embedded in the resonance structure, exhibiting hierarchical organization and internal bifurcation patterns. These features highlight the role of critical orbit interactions in shaping both global and local aspects of the parameter space.

In the large-degree regime, the global shape of the connectedness loci appears to stabilize and approach that of cubic polynomial slices, while increasingly fine structures accumulate near singular parameters. This suggests a form of geometric universality, although a precise theoretical description remains an interesting open problem.

Taken together, these observations indicate that the parameter spaces considered here are governed by a subtle interplay between critical orbit relations, bifurcation mechanisms, and arithmetic effects. They point to deep connections with classical polynomial dynamics, while also revealing new phenomena specific to the rational setting, thereby opening several directions for future investigation.

\section*{Acknowledgements}
The author thanks the Institute for Pure and Applied Mathematics (IMPA, Rio de Janeiro) for its support and hospitality during the 2024 Postdoctoral Summer Program, where the first simulations were carried out. The author was partially supported by  CNPq 169201/2023-6 grant.



\FloatBarrier

\bibliographystyle{plain} 
\bibliography{sample.bib}

@article{baranski1998newton,
  title={From Newton’s method to exotic basins Part I: The parameter space},
  author={Bara{\'n}ski, Krzysztof},
  journal={Fundamenta Mathematicae},
  volume={158},
  number={3},
  pages={249--288},
  year={1998}
}

@article{baranski2001newton,
  title={From Newton's method to exotic basins Part II: Bifurcation of the Mandelbrot-like sets},
  author={Bara{\'n}ski, Krzysztof},
  journal={Fundamenta Mathematicae},
  volume={1},
  number={168},
  pages={1--55},
  year={2001}
}

@article{zhang2025parabolic,
  title={Parabolic Implosion in the Parameter Space of Cubic Polynomials},
  author={Zhang, Runze},
  journal={arXiv preprint arXiv:2508.16430},
  year={2025}
}

@article{zhang2024dynamical,
  title={On dynamical parameter space of cubic polynomials with a parabolic fixed point},
  author={Zhang, Runze},
  journal={Journal of the London Mathematical Society},
  volume={110},
  number={6},
  pages={e70038},
  year={2024},
  publisher={Wiley Online Library}
}

@article{petersen2025lemon,
  title={Lemon limbs of the cubic connectedness locus},
  author={Petersen, Carsten Lunde and Zakeri, Saeed},
  journal={arXiv preprint arXiv:2504.19081},
  year={2025}
}

@article{milnor1992remarks,
  title={Remarks on iterated cubic maps},
  author={Milnor, John},
  journal={Experimental Mathematics},
  volume={1},
  number={1},
  pages={5--24},
  year={1992},
  publisher={Taylor \& Francis}
}

@article{KawahiraKisaka2026,
  author  = {Kawahira, Tomoki and Kisaka, Masashi},
  title   = {Julia sets appear quasiconformally in the Mandelbrot set},
  journal = {Conformal Geometry and Dynamics},
  volume  = {30},
  year    = {2026},
  pages   = {1--43},
  doi     = {10.1090/ecgd/401},
  note    = {See also arXiv:1804.00176}
}

@article{coronel2018non,
  title={Non computable Mandelbrot-like sets for a one-parameter complex family},
  author={Coronel, Daniel and Rojas, Cristobal and Yampolsky, Michael},
  journal={Information and Computation},
  volume={262},
  pages={110--122},
  year={2018},
  publisher={Elsevier}
}

@article{buff2013limits,
  title={Limits of Degenerate Parabolic Quadratic Rational Maps},
  author={Buff, Xavier and {\'E}calle, Jean and Epstein, Adam},
  journal={Geometric and Functional Analysis},
  volume={23},
  number={1},
  pages={42--95},
  year={2013},
  publisher={Springer}
}

@article{milnor20099,
  title={9 Cubic Polynomial Maps with Periodic Critical Orbit, Part I},
  author={Milnor, John},
  journal={Complex Dynamics: Families and Friends},
  pages={333--411},
  year={2009},
  publisher={AK Peters, Ltd.}
}

@article{milnor1993geometry,
  title={Geometry and dynamics of quadratic rational maps},
  author={Milnor, John and Milnor and Lei, Tan},
  journal={Experimental mathematics},
  volume={2},
  number={1},
  pages={37--83},
  year={1993},
  publisher={Taylor \& Francis}
}

@article{buff2001julia,
  title={Julia sets in parameter spaces},
  author={Buff, Xavier and Henriksen, Christian},
  journal={Communications in Mathematical Physics},
  volume={220},
  number={2},
  pages={333--375},
  year={2001},
  publisher={Springer}
}

@article{canela2015family,
  title={On a family of rational perturbations of the doubling map},
  author={Canela, Jordi and Fagella, N{\'u}ria and Garijo, Antonio},
  journal={Journal of Difference Equations and Applications},
  volume={21},
  number={8},
  pages={715--741},
  year={2015},
  publisher={Taylor \& Francis}
}

@article{zakeri1999dynamics,
  title={Dynamics of cubic Siegel polynomials},
  author={Zakeri, Saeed},
  journal={Communications in mathematical physics},
  volume={206},
  number={1},
  pages={185--233},
  year={1999},
  publisher={Springer}
}

@article{mcmullen2000mandelbrot,
  title={The Mandelbrot set is universal},
  author={MCMULLEN, C},
  journal={The Mandelbrot Set, Theme and Variations, London Mathematical Society Lecture Note Series},
  volume={274},
  pages={1--18},
  year={2000}
}

@article{douady1987disques,
  title={Disques de Siegel et anneaux de Herman},
  author={Douady, Adrien},
  journal={Ast{\'e}risque},
  volume={152},
  number={153},
  pages={4},
  year={1987}
}

@article{milnor2000rational,
  title={On rational maps with two critical points},
  author={Milnor, John},
  journal={Experimental Mathematics},
  volume={9},
  number={4},
  pages={481--522},
  year={2000},
  publisher={Taylor \& Francis}
}

@article{bonifant2010cubic,
  title={Cubic polynomial maps with periodic critical orbit, part ii: escape regions},
  author={Bonifant, Araceli and Kiwi, Jan and Milnor, John},
  journal={Conformal Geometry and Dynamics of the American Mathematical Society},
  volume={14},
  number={4},
  pages={68--112},
  year={2010}
}

@article{bonifant2025cubic,
  title={Cubic Polynomial Maps with Periodic Critical Orbit, Part III: Tessellations and Orbit Portraits},
  author={Bonifant, Araceli and Milnor, John},
  journal={arXiv preprint arXiv:2503.08868},
  year={2025}
}

@article{DouadyHubbard,
  author = {Douady, Adrien and Hubbard, John H.},
  title = {On the dynamics of polynomial-like mappings},
  journal = {Annales Scientifiques de l'École Normale Supérieure},
  volume = {18},
  year = {1985},
  pages = {287--343}
}

@article{Milnor2000,
  author = {Milnor, John},
  title = {Geometry and dynamics of quadratic rational maps},
  journal = {Experimental Mathematics},
  volume = {2},
  year = {1993},
  pages = {37--83}
}

@article{lomonaco2015parabolic,
  title={Parabolic-like mappings},
  author={Lomonaco, Luna},
  journal={Ergodic theory and dynamical systems},
  volume={35},
  number={7},
  pages={2171--2197},
  year={2015},
  publisher={Cambridge University Press}
}

@article{inou2012combinatorics,
  title={Combinatorics and topology of straightening maps, I: Compactness and bijectivity},
  author={Inou, Hiroyuki and Kiwi, Jan},
  journal={Advances in Mathematics},
  volume={231},
  number={5},
  pages={2666--2733},
  year={2012},
  publisher={Elsevier}
}

@article{inouvisualization,
  title={VISUALIZATION OF THE BIFURCATION LOCUS OF CUBIC POLYNOMIAL FAMILY},
  author={INOU, HIROYUKI}
}

@article{bonifant2023relations,
  title={Relations Between Escape Regions in the Parameter Space of Cubic Polynomials},
  author={Bonifant, Araceli and Estabrooks, Chad and Sharland, Thomas},
  journal={Arnold Mathematical Journal},
  volume={9},
  number={2},
  pages={245--265},
  year={2023},
  publisher={Springer}
}

@article{gallas1993structure,
  title={Structure of the parameter space of the H{\'e}non map},
  author={Gallas, Jason AC},
  journal={Physical review letters},
  volume={70},
  number={18},
  pages={2714},
  year={1993},
  publisher={APS}
}

@article{gallas2010structure,
  title={The structure of infinite periodic and chaotic hub cascades in phase diagrams of simple autonomous flows},
  author={Gallas, Jason AC},
  journal={International Journal of Bifurcation and Chaos},
  volume={20},
  number={02},
  pages={197--211},
  year={2010},
  publisher={World Scientific}
}

\vspace{1.5em}

\noindent Instituto de Matemática Pura e Aplicada\\
Estrada Dona Castorina 110, Jardim Botânico\\
22460-320, Rio de Janeiro, Brazil\\
\texttt{ivan.suarez@impa.br}


\end{document}